# OPTIMAL POINTWISE APPROXIMATION OF SDES BASED ON BROWNIAN MOTION AT DISCRETE POINTS

By Thomas Müller-Gronbach

*Otto-von-Guericke-Universität Magdeburg*

We study pathwise approximation of scalar stochastic differential equations at a single point. We provide the exact rate of convergence of the minimal errors that can be achieved by arbitrary numerical methods that are based (in a measurable way) on a finite number of sequential observations of the driving Brownian motion. The resulting lower error bounds hold in particular for all methods that are implementable on a computer and use a random number generator to simulate the driving Brownian motion at finitely many points. Our analysis shows that approximation at a single point is strongly connected to an integration problem for the driving Brownian motion with a random weight. Exploiting general ideas from estimation of weighted integrals of stochastic processes, we introduce an adaptive scheme, which is easy to implement and performs asymptotically optimally.

**1. Introduction.** We consider a scalar stochastic differential equation

$$dX(t) = a(t, X(t))\,dt + \sigma(t, X(t))\,dW(t), \qquad t \in [0,1], \tag{1}$$

with initial value $X(0)$. Here $W$ denotes a one-dimensional Brownian motion, and $a:[0,1]\times\mathbb{R}\to\mathbb{R}$ and $\sigma:[0,1]\times\mathbb{R}\to\mathbb{R}$ satisfy standard smoothness conditions.

In most cases an explicit solution of (1) will not be available so that an approximation $\widehat{X}$ must be used. Assume that the driving Brownian motion $W$ may be evaluated at a finite number of points. Then the following questions are of interest:

1. Where in the unit interval should these evaluations be made and how should the resulting data be used in order to obtain the best possible approximation to the solution?

2. What is the minimal error that can be achieved if at most $N$ evaluations of $W$ are made on the average?

The analysis of these problems clearly needs the specification of an error criterion. The two main approaches in the literature are:

(i) *approximation at a finite number of points*, that is, $\widehat{X}$ is compared to the solution $X$ at finitely many points in the unit interval,

(ii) *global approximation*, that is, $\widehat{X}$ is compared to the solution $X$ globally on the unit interval.

First results for global approximation are due to Pardoux and Talay (1985) who studied almost surely uniform convergence of specific approximations. Faure (1992) determines an upper bound with an unspecified constant for the average $L_\infty$-error of a Euler scheme with piecewise linear interpolation. Complete answers (in an asymptotic sense) to the questions 1 and 2 above are given in Hofmann, Müller-Gronbach and Ritter (2001) for the average $L_2$-error and Müller-Gronbach (2002b) for the average $L_\infty$-error. In these papers the exact rate of convergence of the minimum error is determined and adaptive methods are presented that are easy to implement and perform asymptotically optimally.

Much less is known for the problem of approximation at a finite number of points. Here, the majority of results deal only with upper bounds for the error of specific schemes at the discretization points; see, for example, Kloeden and Platen (1995) for an overview. Lower bounds for approximation at $t=1$ were first presented in Clark and Cameron (1980) who considered an autonomous equation (1) with constant diffusion $\sigma = 1$ and determined the rate of convergence of the minimal mean squared error that can be obtained by equidistant evaluation of the driving Brownian motion $W$. Rümelin (1982) studied autonomous equations with a nonconstant diffusion coefficient and presented the order of the minimal error that can be obtained by Runge–Kutta methods based on equidistant evaluation of $W$. The most fargoing result is due to Cambanis and Hu (1996) who analyzed the mean squared error of the conditional expectation of $X(1)$ given observations of $W$ at points that are regularly generated by some density. They provided the rate of convergence of the corresponding mean squared error and determined the optimal density. Clearly, all these results only provide partial answers to the above questions 1 and 2. For instance, the implementation of a conditional expectation will be a hard task in general. Furthermore, considerations are restricted to numerical methods that are based on sampling $W$ at prefixed points in the unit interval (either equidistant or regularly generated by some density). Adaptive methods which take into account the particular trajectory of the solution are not covered. See Remarks 1 and 3 for a discussion.

In the present paper we provide a detailed analysis of approximation at $t=1$ with respect to the questions 1 and 2. Our results cover all numerical methods that are based on the initial value $X(0)$ and finitely many sequential observations

$$W(\tau_1),\ldots,W(\tau_\nu)$$

of the driving Brownian motion $W$. Except for measurability conditions, we do not impose any further restrictions. The $k$th evaluation point $\tau_k$ may depend on the previous evaluations $X(0),W(\tau_1),\ldots,W(\tau_{k-1})$ and the total number $\nu$ of observations of $W$ may be determined by a stopping rule. Finally, the resulting discrete data may be used in any way to generate an estimator

$$\widehat{X}(1)=\phi(W(\tau_1),\ldots,W(\tau_\nu))$$

of $X(1)$, the solution at $t=1$. For example, the adaptive Euler–Maruyama scheme recently introduced in Lamba, Mattingly and Stuart (2003) is of this type.

The error of $\widehat{X}(1)$ is defined by

$$e_p(\widehat{X}(1))=(E|X(1)-\widehat{X}(1)|^p)^{1/p},$$

where $p\in[1,\infty[$, and $c(\widehat{X}(1))=E(\nu)$ is the average number of evaluations of $W$ used by $\widehat{X}(1)$.

Our analysis shows that the problem of pathwise approximation at $t=1$ is strongly connected to an integration problem for the driving Brownian motion $W$ with the random weight

$$\mathcal{Y}(t)=\mathcal{M}(t,1)\cdot\mathcal{G}(t,X(t)),\qquad t\in[0,1],$$

where $\mathcal{G}=\sigma a^{(0,1)}-\sigma^{(1,0)}-a\sigma^{(0,1)}-1/2\cdot\sigma^2\sigma^{(0,2)}$ involves partial derivatives of $a$ and $\sigma$, and the one-dimensional random field $\mathcal{M}$ is given by

$$\mathcal{M}(t,s)=\exp\biggl(\int_t^s(a^{(0,1)}-1/2\cdot(\sigma^{(0,1)})^2)(u,X(u))\,du + \int_t^s\sigma^{(0,1)}(u,X(u))\,dW(u)\biggr)$$

for $0\le t\le s\le 1$. Roughly speaking, $\mathcal{M}(t,\cdot)$ is the $L_2$-derivative of the solution $X$ with respect to its state at time $t$; see Remark 6.

To give a flavor of our results, let $p=2$, and consider the minimal error

$$e_2^{**}(N)=\inf\{e_2(\widehat{X}(1)):c(\widehat{X}(1))\le N\}$$

that can be achieved by numerical methods using at most $N$ evaluations of the driving Brownian motion on the average. By Theorem 1(i),

$$\lim_{N\to\infty}N\cdot e_2^{**}(N)=1/\sqrt{12}\cdot\biggl(E\biggl(\int_0^1|\mathcal{Y}(t)|^{2/3}\,dt\biggr)\biggr)^{3/2}, \tag{2}$$

which answers question 2 in an asymptotic sense.

For answering question 1 we exploit general ideas from estimation of weighted integrals of stochastic processes; see, for example, Ritter (2000) and the references therein. We construct an easy to implement adaptive scheme $\widehat{X}^{**}_{2,n}$ with step-size roughly proportional to $|\widehat{\mathcal{Y}}_n(t)|^{-2/3}$, where $\widehat{\mathcal{Y}}_n$ is a suitable approximation to the random weight $\mathcal{Y}$. The resulting approximation $\widehat{X}^{**}_{2,n}(1)$ at $t=1$ satisfies

$$\lim_{n\to\infty} c(\widehat{X}^{**}_{2,n}(1)) \cdot e_2(\widehat{X}^{**}_{2,n}(1)) = 1/\sqrt{12} \cdot \left(E\left(\int_0^1 |\mathcal{Y}(t)|^{2/3}\, dt\right)\right)^{3/2}; \tag{3}$$

see Theorem 2(i). Consequently, by (2) this method performs asymptotically optimally for every equation (1) with a nonzero asymptotic constant on the right-hand side above.

A natural question is whether the asymptotic constant in (3) can also be achieved by a numerical method based on a prefixed discretization. The answer turns out to be negative in general. In fact, consider the minimal error

$$e_2(N) = \inf\{e_2(E(X(1)|W(t_1),\dots,W(t_N))) : 0 \le t_1 \le \cdots \le t_N \le 1\}$$

that can be obtained if the driving Brownian motion $W$ is evaluated at $N$ prefixed points in the unit interval. By Theorem 1(iii),

$$\lim_{N\to\infty} N \cdot e_2(N) = 1/\sqrt{12} \cdot \left(\int_0^1 (E|\mathcal{Y}(t)|^2)^{1/3}\, dt\right)^{3/2}. \tag{4}$$

Thus the order of convergence is still $1/N$ but the asymptotic constant in (4) may be considerably larger than the asymptotic constant in (2); see Example 1. Somewhat surprisingly, as a by-product of (4), it turns out that in general the Milstein scheme does not perform asymptotically optimally; see Remark 7.

In Section 2 we state our assumptions on equation (1). We use global Lipschitz and linear growth conditions on the drift coefficient $a$, the diffusion coefficient $\sigma$ and partial derivatives of these coefficients, as well as a moment condition on the initial value $X(0)$.

Best rates of convergence for approximation at $t=1$ based on point evaluations of $W$ are stated in Section 3. More specifically, we analyze the minimal errors that can be achieved if $W$ is evaluated at:

(a) sequentially chosen points $\tau_1,\dots,\tau_\nu$ with $E(\nu) \le N$,
(b) sequentially chosen points $\tau_1,\dots,\tau_\nu$ with $\nu \le N$,
(c) prefixed points $t_1,\dots,t_N$,
(d) equidistant points $1/N, 2/N, \dots, 1$.

In Section 4 we introduce a new class of numerical schemes, which leads to asymptotically optimal approximations for each of the cases (a)–(d) above. Proofs are postponed to Section 5 and the Appendix.

**2. Assumptions.** We will use the following Lipschitz and linear growth conditions on functions $f:[0,1]\times\mathbb{R}\to\mathbb{R}$.

(L) There exists a constant $K>0$ such that

$$|f(t,x)-f(t,y)|\le K\cdot|x-y|$$

for all $t\in[0,1]$ and $x,y\in\mathbb{R}$.

(LG) There exists a constant $K>0$ such that

$$|f(t,x)|\le K\cdot(1+|x|)$$

for all $t\in[0,1]$ and $x\in\mathbb{R}$.

(LLG) There exists a constant $K>0$ such that

$$|f(s,x)-f(t,x)|\le K\cdot(1+|x|)\cdot|s-t|$$

for all $s,t\in[0,1]$ and $x\in\mathbb{R}$.

Throughout this paper we impose the following regularity conditions on the drift coefficient $a$, the diffusion coefficient $\sigma$ and the initial value $X(0)$.

(A) (i) Both $a$ and $\sigma$ satisfy (L) as well as (LLG).
(ii) The partial derivatives

$$a^{(1,0)},\quad a^{(0,1)},\quad a^{(0,2)},\quad \sigma^{(1,0)},\quad \sigma^{(0,1)},\quad \sigma^{(0,2)}$$

exist and satisfy (L) as well as (LLG).
(iii) The functions $\sigma^2a^{(0,2)}$ and $\sigma^2\sigma^{(0,2)}$ satisfy (LG).
(iv) The function $\sigma\sigma^{(0,2)}$ satisfies (L).

(B) The initial value $X(0)$ is independent of $W$ and satisfies $E|X(0)|^{16p}<\infty$.

For instance, (A) is satisfied if the partial derivatives

$$a^{(i,j)},\quad \sigma^{(i,j)},\qquad i=0,1,2,\ j=0,1,2,3,$$

exist and are continuous and bounded.

Note that (A) together with (B) implies that a pathwise unique strong solution of the equation (1) with initial value $X(0)$ exists. In particular, the conditions assure that

$$E\Big(\sup_{0\le t\le 1}|X(t)|^{16p}\Big)<\infty \tag{5}$$

as well as

$$E|X(s)-X(t)|^{16p}\le c\cdot|s-t|^{8p}, \tag{6}$$

where the constant $c>0$ only depends on $p$ and the constants from (A) and (B).

**3. Best rates of convergence.** We consider arbitrary numerical methods for pathwise approximation of the solution $X$ at the point $t=1$ that are based on a realization of the initial value $X(0)$ and a finite number of observations of a trajectory of the driving Brownian motion $W$ at points in the unit interval. The formal definition of the class of these methods and subclasses of interest is given in Section 3.1. Section 3.2 contains the analysis of the corresponding minimal errors.

3.1. *General methods for approximation at* $t=1$. A general adaptive approximation $\widehat{X}(1)$ of $X(1)$ is defined by three sequences

$$\psi = (\psi_n)_{n\in\mathbb{N}}, \qquad \chi = (\chi_n)_{n\in\mathbb{N}}, \qquad \phi = (\phi_n)_{n\in\mathbb{N}},$$

of measurable mappings

$$\begin{aligned}\psi_n &: \mathbb{R}^n \to ]0,1],\\ \chi_n &: \mathbb{R}^{n+1} \to \{\mathrm{STOP}, \mathrm{GO}\},\\ \phi_n &: \mathbb{R}^{n+1} \to \mathbb{R}.\end{aligned}$$

The sequence $\psi$ determines the evaluation sites of a trajectory $w$ of $W$ in the interval $]0,1]$. The total number of evaluations to be made is determined by the sequence $\chi$ of stopping rules. Finally, $\phi$ is used to obtain the real-valued approximation to the solution $X$ at $t=1$ from the observed data.

To be more precise, the sequential observation of a trajectory $w$ starts at the knot $\psi_1(x)$, where $x$ denotes the realization of the initial value. After $n$ steps we have obtained the data $D_n(x,w) = (x, y_1, \dots, y_n)$, where $y_1 = w(\psi_1(x)), \dots, y_n = w(\psi_n(x, y_1, \dots, y_{n-1}))$, and we decide to stop or to further evaluate $w$ according to the value of $\chi_n(D_n(x,w))$. The total number of observations is thus given by $\nu(x,w) = \min\{n \in \mathbb{N} : \chi_n(D_n(x,w)) = \mathrm{STOP}\}$. If $\nu(x,w) < \infty$, then the whole data $D(x,w) = D_{\nu(x,w)}(x,w)$ are used to construct the estimate $\phi_{\nu(x,w)}(D(x,w)) \in \mathbb{R}$.

For obvious reasons we require $\nu(X(0), W) < \infty$ with probability 1. Then the resulting approximation is given by

$$\widehat{X}(1) = \phi_{\nu(X(0),W)}(D(X(0),W)).$$

As a rough measure for the computational cost of $\widehat{X}(1)$ we use

$$c(\widehat{X}(1)) = E(\nu(X(0),W)),$$

that is, the expected number of evaluations of the driving Brownian motion $W$. Clearly, a more realistic measure also involves, for example, a count of the arithmetical operations needed to compute $\widehat{X}(1)$.

Let $\mathcal{X}^{**}$ denote the class of all methods of the above form and put

$$\mathcal{X}^{**}_N = \{\widehat{X}(1) \in \mathcal{X}^{**} : c(\widehat{X}(1)) \le N\}.$$

Then

$$e_p^{**}(N) = \inf\{e_p(\widehat{X}(1)) : \widehat{X}(1) \in \mathcal{X}_N^{**}\}$$

is the minimal error that can be obtained by approximations that use at most $N$ sequential observations of $W$ on the average.

The number and the location of the evaluation sites that are used by an approximation $\widehat{X}(1) \in \mathcal{X}^{**}$ depend on the respective realization $x$ of the initial value $X(0)$ and the path $w$ of the driving Brownian motion $W$. It is natural to ask whether, in general, the minimal errors $e_p^{**}(N)$ can (asymptotically) be achieved by methods that use the same evaluation sites for every trajectory of $W$. In order to investigate questions of this type, we introduce the following subclasses of $\mathcal{X}^{**}$ that are subject to certain restrictions on the choice of evaluation sites.

The subclass $\mathcal{X}^* \subset \mathcal{X}^{**}$ consists of all approximations that use the same number of observations for every $x$ and $w$. Formally, this means that the mappings $\chi_n$ are constant and $\nu = \min\{n \in \mathbb{N} : \chi_n = \text{STOP}\}$.

Additionally, we consider the subclass $\mathcal{X} \subset \mathcal{X}^*$ of all approximations that evaluate $W$ at the same points for every $x$ and every path $w$. Formally, the mappings $\psi_n$ and $\chi_n$ are constant so that $\nu = \min\{n \in \mathbb{N} : \chi_n = \text{STOP}\}$ and $D(x, w) = (x, w(\psi_1), \ldots, w(\psi_\nu))$. For instance, if the discretization is fixed, then the corresponding Euler scheme and the Milstein scheme at $t = 1$ belong to the class $\mathcal{X}$.

Finally, the class $\mathcal{X}^{\text{equi}} \subset \mathcal{X}$ consists of all approximations that use equidistant evaluation sites for the driving Brownian motion.

The definition of the respective classes $\mathcal{X}_N^*$, $\mathcal{X}_N$, $\mathcal{X}_N^{\text{equi}}$ and the corresponding minimal errors $e_p^*(N)$, $e_p(N)$ and $e_p^{\text{equi}}(N)$ is canonical.

We stress that the class $\mathcal{X}^{**}$ contains all commonly studied methods for approximation at $t = 1$ that are based on function values of the driving Brownian motion. Formally, the corresponding sequences $\psi$, $\chi$ and $\phi$ then depend on the respective drift coefficient $a$ and diffusion coefficient $\sigma$. In the majority of cases, partial information about the coefficients, for example, finitely many function values or derivative values, are sufficient to compute the approximations $\widehat{X}(1)$.

In the present paper we present (asymptotically) sharp upper and lower bounds for the minimal errors defined above. The upper bounds are achieved by methods that also need only partial information about $a$ and $\sigma$. On the other hand, no restriction on the available information about $a$ and $\sigma$ is present in the definition of the class $\mathcal{X}^{**}$. Therefore, the lower bounds hold even for strong approximations that may specifically be tuned to the respective coefficients. As an example, consider an approximation of the form

$$\widehat{X}(1) = E(X(1) | W(t_1), \ldots, W(t_N)),$$

which belongs to the class $\mathcal{X}_N$ and might even not be implementable.

3.2. *Analysis of minimal errors.* Let $\mathfrak{m}_p$ denote the $p$th root of the $p$th absolute moment of a standard normal variable, that is,

$$\mathfrak{m}_p = \left( \int_{-\infty}^{\infty} |y|^p/(2\pi)^{1/2} \cdot \exp(-y^2/2)\, dy \right)^{1/p}.$$

Recall the weighting process $\mathcal{Y}$ from Section 1 and define the constants

$$C_p^{**} = \mathfrak{m}_p \cdot \left( E\left( \int_0^1 |\mathcal{Y}(t)|^{2/3}\, dt \right)^{3p/2(p+1)} \right)^{(p+1)/p},$$

$$C_p^{*} = \mathfrak{m}_p \cdot \left( E\left( \int_0^1 |\mathcal{Y}(t)|^{2/3}\, dt \right)^{3p/2} \right)^{1/p},$$

$$C_2 = \left( \int_0^1 (E|\mathcal{Y}(t)|^2)^{1/3}\, dt \right)^{3/2},$$

$$C_p^{\mathrm{equi}} = \mathfrak{m}_p \cdot \left( E\left( \int_0^1 |\mathcal{Y}(t)|^2\, dt \right)^{p/2} \right)^{1/p}.$$

Theorem 1. *The minimal errors satisfy:*

(i) $\lim_{N\to\infty} N \cdot e_p^{**}(N) = C_p^{**}/\sqrt{12}$,
(ii) $\lim_{N\to\infty} N \cdot e_p^{*}(N) = C_p^{*}/\sqrt{12}$,
(iii) $\lim_{N\to\infty} N \cdot e_2(N) = C_2/\sqrt{12}$,
(iv) $\lim_{N\to\infty} N \cdot e_p^{\mathrm{equi}}(N) = C_p^{\mathrm{equi}}/\sqrt{12}$.

Clearly, the asymptotic constants vanish altogether iff $C_2^{\mathrm{equi}} = 0$. Thus, if $C_2^{\mathrm{equi}} > 0$, then the order of convergence of the minimal errors is $1/N$ for all of the above classes. However, note that

$$C_p^{**} \le C_p^{*} \le C_p^{\mathrm{equi}}, \qquad C_2^{**} \le C_2^{*} \le C_2 \le C_2^{\mathrm{equi}},$$

with strict inequality in most cases. See Remark 2 for the case $C_2^{\mathrm{equi}} = 0$ and Remark 4 for a characterization of equality of the asymptotic constants.

Example 1. Consider the linear equation

$$dX(t) = \alpha(t) \cdot X(t)\, dt + \beta(t) \cdot X(t)\, dW(t)$$

with initial condition $X(0) = 1$. Clearly, condition (A) is satisfied if $\alpha$ and $\beta$ have Lipschitz continuous derivatives $\alpha'$ and $\beta'$, respectively. The corresponding field $\mathcal{M}$ is given by

$$\mathcal{M}(t, s) = \exp\left( \int_t^s (\alpha - 1/2 \cdot \beta^2)(u)\, du + \int_t^s \beta(u)\, dW(u) \right),$$

and we have $X(t) = \mathcal{M}(0,t)$. Moreover, $\mathcal{G}(t,x) = -\beta'(t)\cdot x$, so that the weighting process $\mathcal{Y}$ satisfies $\mathcal{Y}(t) = -\beta'(t)\cdot\mathcal{M}(0,t)\cdot\mathcal{M}(t,1) = -\beta'(t)\cdot X(1)$. Straightforward calculations yield for $q\in\mathbb{R}\setminus\{0\}$,

$$(E|X(1)|^q)^{1/q} = e^{\|\alpha\|_1 - 1/2\cdot\|\beta\|_2^2}\cdot e^{q/2\cdot\|\beta\|_2^2}.$$

Thus

$$\begin{aligned}
C_p^{**} &= \mathfrak{m}_p\cdot e^{\|\alpha\|_1 - 1/2\cdot\|\beta\|_2^2}\cdot\|\beta'\|_{2/3}\cdot e^{p/2(p+1)\cdot\|\beta\|_2^2},\\
C_p^{*} &= \mathfrak{m}_p\cdot e^{\|\alpha\|_1 - 1/2\cdot\|\beta\|_2^2}\cdot\|\beta'\|_{2/3}\cdot e^{p/2\cdot\|\beta\|_2^2},\\
C_2 &= C_2^*,\\
C_p^{\text{equi}} &= \mathfrak{m}_p\cdot e^{\|\alpha\|_1 - 1/2\cdot\|\beta\|_2^2}\cdot\|\beta'\|_2\cdot e^{p/2\cdot\|\beta\|_2^2}.
\end{aligned}$$

If $\alpha = 0$ and $\beta(t) = b\cdot t$ with $b\in\mathbb{R}$, then

$$C_2^{\text{equi}} = C_2 = C_2^* = |b|\cdot e^{b^2/6}$$

and

$$C_2^{**} = |b|\cdot e^{-2b^2/9},$$

which shows that adapting the number of evaluations of $W$ to the particular trajectory of the solution $X$ is essential in this case. Note that the constant $C_2^{**}$ is achieved by the adaptive method to be introduced in Section 4.3.1. Thus, if, for example, $b=5$, then, asymptotically, the error of this method is at least 1/258 times smaller than the error of any approximation based on a fixed number of evaluations of $W$.

Remark 1. Clark and Cameron (1980) consider the autonomous equation

$$dX(t) = a(X(t))\,dt + dW(t), \qquad X(0) = x\in\mathbb{R},$$

where $a$ has bounded derivatives up to order 3. They obtain

$$\begin{aligned}
&\lim_{N\to\infty} N\cdot(E|X(1) - E(X(1)|W(1/N),\dots,W(1))|^2)^{1/2}\\
&\quad = \left(\int_0^1 E\left((a'(X(t)))^2\cdot\exp\left(2\cdot\int_t^1 a'(X(s))\,ds\right)\right)dt\right)^{1/2}\Big/\sqrt{12}.
\end{aligned}$$

Note that the corresponding weighting process is given by

$$\mathcal{Y}(t) = a'(X(t))\cdot\exp\left(\int_t^1 a'(X(s))\,ds\right)$$

so that the above result is a consequence of Theorem 1(iv).

More generally, Cambanis and Hu (1996) study autonomous equations

$$dX(t) = a(X(t))\,dt + \sigma(X(t))\,dW(t), \qquad X(0) = x \in \mathbb{R},$$

where $a$ and $\sigma$ have bounded derivatives up to order 3. They analyze the minimum error that can be achieved by methods from the class $\mathcal{X}$ that are based on so-called regularly generated discretizations.

To be more precise, let $h$ be a strictly positive density on $[0,1]$ and define the discretization

$$0 < t_1^{(h)} < \cdots < t_N^{(h)} = 1$$

by taking the $l/N$-quantiles corresponding to $h$, that is,

$$\int_0^{t_l^{(h)}} h(t)\,dt = l/N, \qquad l = 1, \dots, N.$$

Consider the optimal approximation in the mean squared sense

$$\widehat{X}_N^{(h)}(1) = E(X(1)|W(t_1^{(h)}), \dots, W(t_N^{(h)}))$$

that is based on the observations $W(t_1^{(h)}), \dots, W(t_N^{(h)})$ and put

$$C(h) = \left(\int_0^1 E|\mathcal{Y}(t)|^2/h^2(t)\,dt\right)^{1/2}.$$

If $h$ has a bounded derivative, then

$$\lim_{N\to\infty} N \cdot e_2(\widehat{X}_N^{(h)}(1)) = C(h)/\sqrt{12}.$$

Taking $h = 1$ yields Theorem 1(iv) in the case $p = 2$, since

$$e_2(\widehat{X}_N^{(1)}(1)) = e_2^{\mathrm{equi}}(N) \quad \text{and} \quad C(1) = C_2^{\mathrm{equi}}.$$

Taking

$$h^*(t) = (E|\mathcal{Y}(t)|^2)^{1/3} \Big/ \int_0^1 (E|\mathcal{Y}(s)|^2)^{1/3}\,ds$$

yields the minimal constant

$$C_2 = C(h^*) = \inf_h C(h).$$

Thus, by Theorem 1(iii), the approximation $\widehat{X}_N^{(h^*)}(1)$ is asymptotically optimal in the class $\mathcal{X}$ if $C_2 > 0$. However, note that the method $\widehat{X}_N^{(h^*)}(1)$ is much harder to implement than the asymptotically optimal method introduced in Section 4.

REMARK 2. Theorem 1 determines the rates of convergence of the minimal errors only in the case of nonzero asymptotic constants. Clearly, these constants vanish altogether iff with probability 1,

$$\mathcal{G}(t, X(t)) = 0 \qquad \text{for every } t \in [0,1]. \tag{7}$$

For a large class of equations, it turns out that (7) holds iff there exists a measurable function $f : \mathbb{R} \times [0,1] \times \mathbb{R} \to \mathbb{R}$ such that, with probability 1,

$$X(t) = f(X(0), t, W(t)) \qquad \text{for every } t \in [0,1]. \tag{8}$$

Obviously, if (8) holds, then $X(1)$ can be simulated exactly. Thus (8) implies (7) by Theorem 1.

Clark and Cameron (1980) provide sufficient conditions for the equivalence of (7) and (8) in the case of autonomous equations. Slightly modifying their approach, one can also treat general equations. If, additionally to assumption (A), the conditions

(i) $a$ and $\sigma^{(1,0)}$ are bounded,
(ii) $\inf_{t,x} |\sigma(t,x)| > 0$,

are satisfied, then (7) and (8) are equivalent.

The equivalence of (7) and (8) also holds for the linear equation from Example 1. Note that condition (ii) above must not be satisfied in this case. However, (7) implies $\beta' = 0$ and therefore

$$X(t) = \exp\left(\int_0^t (\alpha - 1/2 \cdot \beta^2)(u)\, du + \beta(t) \cdot W(t)\right), \qquad t \in [0,1].$$

Finally, assume that $a$ and $\sigma$ have partial derivatives of any order. Then, by a general result of Yamato (1979), (8) is equivalent to

$$\mathcal{G} = 0, \tag{6*}$$

which clearly implies (7).

Note that (6*) implies that the Wagner–Platen scheme only uses function values of the driving Brownian motion; see Section 4. Thus the order of convergence of the minimal errors $e_p^{\mathrm{equi}}(N)$ is at least $1/N^{3/2}$ in this case.

REMARK 3. Rümelin (1982) analyzes a class $\overline{\mathcal{X}}$ of Runge–Kutta methods based on an equidistant discretization, that is,

$$\overline{\mathcal{X}} \subset \mathcal{X}^{\mathrm{equi}},$$

with respect to the mean squared error at $t = 1$, that is, $p = 2$. For this class Rümelin shows that, under stronger conditions on $a$ and $\sigma$, the order of convergence of the corresponding minimal errors is $1/N$ iff $\mathcal{G} \neq 0$. Moreover, if $\mathcal{G} = 0$, then the order is at least $1/N^{3/2}$.

REMARK 4. We briefly comment on equality of the asymptotic constants in the case $p=2$. Clearly, $C_2 = C_2^{\text{equi}}$ iff there exists $\gamma \in \mathbb{R}$ such that

$$E(\mathcal{Y}(t))^2 = \gamma \qquad \text{for all } t \in [0,1].$$

Furthermore, $C_2^* = C_2$ iff there exist $t_0 \in [0,1]$ and a function $\gamma \in C([0,1])$ such that, with probability 1,

$$\mathcal{Y}(t) = \gamma(t) \cdot \mathcal{Y}(t_0) \qquad \text{for all } t \in [0,1].$$

Note that the latter condition holds for the linear equation from Example 1 with $\gamma = -\beta'/\beta'(1)$ and $t_0 = 1$.

Finally, by the Markov property of $X$, we have $C_2^{**} = C_2^*$ iff there exists a function $\gamma \in C([0,1])$ such that, with probability 1,

$$\mathcal{Y}(t) = \gamma(t) \qquad \text{for all } t \in [0,1].$$

In particular, if $a$ and $\sigma$ are state independent, then $C_2^{**} = C_2^* = C_2 = \|\sigma'\|_{2/3}$.

REMARK 5. Theorem 1 shows that pathwise approximation at a single point is strongly connected to weighted integration of a Brownian motion. To be more precise, let $\rho\colon [0,1] \to [0,\infty[$ be continuous, and consider the problem of estimating the weighted integral

$$I = \int_0^1 \rho(t) \cdot \widetilde{W}(t)\, dt$$

of a Brownian motion $\widetilde{W}$ on the basis of $N$ observations of $\widetilde{W}$ in the unit interval. The corresponding minimum mean squared error

$$\varepsilon(N) = \inf\{(E(I - E(I|\widetilde{W}(t_1), \ldots, \widetilde{W}(t_N)))^2)^{1/2} : 0 \le t_1 \le \cdots \le t_N \le 1\}$$

satisfies

$$\lim_{N \to \infty} N \cdot \varepsilon(N) = 1/\sqrt{12} \cdot c_\rho,$$

where

$$c_\rho = \left( \int_0^1 (\rho(t))^{2/3}\, dt \right)^{3/2};$$

see Ritter (2000) and the references therein.

Taking the weight

$$\rho(t) = (E(\mathcal{Y}^2(t)))^{1/2}, \qquad t \in [0,1],$$

yields the constant $C_2$ in Theorem 1(iii).

Using the random weight $|\mathcal{Y}|$, we obtain the random constant $c_{|\mathcal{Y}|}$, and

$$C_2^{**} = (E(c_{|\mathcal{Y}|}^{2/3}))^{3/2}.$$

As an illustrating example, consider the linear equation with additive noise

$$dX(t) = a(t)\,dt + \sigma(t)\,dW(t).$$

Then

$$X(1) = X(0) + \int_0^1 a(t)\,dt + \sigma(1)\cdot W(1) - \int_0^1 \sigma'(t)\cdot W(t)\,dt.$$

Since $X(0)$ and $\sigma(1)\cdot W(1)$ can be observed, we are basically dealing with the approximation of the last integral on the right-hand side above. Clearly, in this case the weighting process $\mathcal{Y}$ is nonrandom with $\mathcal{Y} = -\sigma'$.

REMARK 6. Consider, for every $x \in \mathbb{R}$ and $t \in [0,1]$, the solution $X_{t,x}$ of the equation

$$dX_{t,x}(s) = a(s, X_{t,x}(s))\,ds + \sigma(s, X_{t,x}(s))\,dW(s), \qquad t \le s \le 1,$$

with initial value $X_{t,x}(t) = x$. As a well-known fact, the distribution of the process $X_{t,x}$ on $C([t,1])$ coincides with the conditional distribution of the solution $X(s)$, $t \le s \le 1$, given $X(t) = x$. Due to condition (A), for every $s \ge t$, there exists the $L_2$-derivative $X'_{t,x}(s)$ of $X_{t,x}(s)$ with respect to the initial value $x$, that is,

$$\lim_{h\to 0} E(1/h\cdot(X_{t,x+h}(s) - X_{t,x}(s)) - X'_{t,x}(s))^2 = 0.$$

Moreover, the process $X'_{t,x}$ is the unique solution of the equation

$$dX'_{t,x}(s) = a^{(0,1)}(s, X_{t,x}(s))\cdot X'_{t,x}(s)\,ds + \sigma^{(0,1)}(s, X_{t,x}(s))\cdot X'_{t,x}(s)\,dW(s), \qquad t \le s \le 1,$$

with initial value $X'_{t,x}(t) = 1$, and is explicitly given by

$$X'_{t,x}(s) = \exp\biggl(\int_t^s (a^{(0,1)} - 1/2\cdot(\sigma^{(0,1)})^2)(u, X_{t,x}(u))\,du + \int_t^s \sigma^{(0,1)}(u, X_{t,x}(u))\,dW(u)\biggr);$$

see, for example, Friedman (1975) and Karatzas and Shreve (1991). Replacing $X_{t,x}$ by the solution $X$ yields the defining equation for the field $\mathcal{M}$.

**4. Asymptotically optimal adaptive schemes.** Let $k \in \mathbb{N}$ and consider the equidistant discretization

$$t_l = l/k, \qquad l = 0, \dots, k. \tag{9}$$

Our adaptive method basically works as follows. First, we evaluate the driving Brownian motion at a coarse grid (9), and we compute a corresponding truncated Wagner–Platen scheme as well as a discrete approximation to the weighting process $\mathcal{Y}$. Following the main idea for nonrandom weighted integration, the latter estimate determines the number and the location of the additional evaluation sites for the driving Brownian motion. The resulting observations are then used to obtain a suitable approximation of the difference between the Wagner–Platen scheme and its truncated version. Finally, we update the truncated Wagner–Platen scheme by adding this approximation.

For convenience we briefly recall the definition of the Wagner–Platen scheme $\widehat{X}_k^{\mathrm{WP}}$ corresponding to the discretization (9). This scheme is defined by $\widehat{X}_k^{\mathrm{WP}}(0) = X(0)$ and

$$\begin{aligned}
\widehat{X}_k^{\mathrm{WP}}(t_{l+1}) = {} & \widehat{X}_k^{\mathrm{WP}}(t_l) + a(t_l, \widehat{X}_k^{\mathrm{WP}}(t_l)) \cdot (t_{l+1} - t_l) \\
& + \sigma(t_l, \widehat{X}_k^{\mathrm{WP}}(t_l)) \cdot (W(t_{l+1}) - W(t_l)) \\
& + 1/2 \cdot (\sigma\sigma^{(0,1)})(t_l, \widehat{X}_k^{\mathrm{WP}}(t_l)) \cdot ((W(t_{l+1}) - W(t_l))^2 - (t_{l+1} - t_l)) \\
& + (\sigma^{(1,0)} + a\sigma^{(0,1)} - 1/2 \cdot \sigma(\sigma^{(0,1)})^2) \\
& \quad \times (t_l, \widehat{X}_k^{\mathrm{WP}}(t_l)) \cdot (W(t_{l+1}) - W(t_l)) \cdot (t_{l+1} - t_l) \\
& + 1/6 \cdot (\sigma(\sigma^{(0,1)})^2 + \sigma^2\sigma^{(0,2)})(t_l, \widehat{X}_k^{\mathrm{WP}}(t_l)) \cdot (W(t_{l+1}) - W(t_l))^3 \\
& + 1/2 \cdot (a^{(1,0)} + aa^{(0,1)} + 1/2 \cdot \sigma^2 a^{(0,2)})(t_l, \widehat{X}_k^{\mathrm{WP}}(t_l)) \cdot (t_{l+1} - t_l)^2 \\
& + \mathcal{G}(t_l, \widehat{X}_k^{\mathrm{WP}}(t_l)) \cdot \int_{t_l}^{t_{l+1}} (W(s) - W(t_l))\, ds
\end{aligned}$$

for $l = 0, \ldots, k-1$; see Wagner and Platen (1978). For the definition of this scheme in the case of a general system of equations, we refer to Kloeden and Platen (1995).

We stress that in general the Wagner–Platen approximation $\widehat{X}_k^{\mathrm{WP}}(1)$ at $t = 1$ does not belong to the class $\mathcal{X}$ since function values as well as integrals of the trajectories of the driving Brownian motion are used.

4.1. *The truncated Wagner–Platen scheme* $\widehat{X}_k^{\mathrm{WPt}}$. Dropping the last summand in the definition of the scheme above, we obtain a truncated version $\widehat{X}_k^{\mathrm{WPt}}$ of the Wagner–Platen scheme that is based, only on function values of the driving Brownian motion. Formally, $\widehat{X}_k^{\mathrm{WPt}}$ is defined by $\widehat{X}_k^{\mathrm{WPt}}(0) = X(0)$ and

$$\begin{aligned}
\widehat{X}_k^{\mathrm{WPt}}(t_{l+1}) = {} & \widehat{X}_k^{\mathrm{WPt}}(t_l) + a(t_l, \widehat{X}_k^{\mathrm{WPt}}(t_l)) \cdot (t_{l+1} - t_l) \\
& + \sigma(t_l, \widehat{X}_k^{\mathrm{WPt}}(t_l)) \cdot (W(t_{l+1}) - W(t_l)) \\
& + 1/2 \cdot (\sigma\sigma^{(0,1)})(t_l, \widehat{X}_k^{\mathrm{WPt}}(t_l)) \cdot ((W(t_{l+1}) - W(t_l))^2 - (t_{l+1} - t_l)) \\
& + (\sigma^{(1,0)} + a\sigma^{(0,1)} - 1/2 \cdot \sigma(\sigma^{(0,1)})^2)
\end{aligned}$$

$$\times (t_l, \widehat{X}_k^{\mathrm{WPt}}(t_l)) \cdot (W(t_{l+1}) - W(t_l)) \cdot (t_{l+1} - t_l)$$
$$+ 1/6 \cdot (\sigma(\sigma^{(0,1)})^2 + \sigma^2 \sigma^{(0,2)})(t_l, \widehat{X}_k^{\mathrm{WPt}}(t_l)) \cdot (W(t_{l+1}) - W(t_l))^3$$
$$+ 1/2 \cdot (a^{(1,0)} + aa^{(0,1)} + 1/2 \cdot \sigma^2 a^{(0,2)})(t_l, \widehat{X}_k^{\mathrm{WPt}}(t_l)) \cdot (t_{l+1} - t_l)^2$$

for $l = 0, \ldots, k-1$.

4.2. *The discrete approximation $\widehat{\mathcal{Y}}_k$ of the random weight $\mathcal{Y}$.* Note that the random field $\mathcal{M}$ satisfies the stochastic differential equations

$$d\mathcal{M}(t,s) = a^{(0,1)}(s, X(s)) \cdot \mathcal{M}(t,s)\, ds + \sigma^{(0,1)}(s, X(s)) \cdot \mathcal{M}(t,s)\, dW(s), \quad t \le s \le 1, \tag{10}$$

with initial value

$$\mathcal{M}(t,t) = 1$$

for every $t \in [0,1]$.

Using the truncated Wagner–Platen estimates, we thus obtain the following Euler-type approximation to the field $\mathcal{M}$. Put

$$\widehat{m}_l = (1 + a^{(0,1)}(t_l, \widehat{X}_k^{\mathrm{WPt}}(t_l)) \cdot (t_{l+1} - t_l) + \sigma^{(0,1)}(t_l, \widehat{X}_k^{\mathrm{WPt}}(t_l)) \cdot (W(t_{l+1}) - W(t_l)))$$

and define the scheme $\widehat{\mathcal{M}}_k$ by

$$\widehat{\mathcal{M}}_k(t_l, t_r) = \begin{cases} \widehat{m}_l \cdots \widehat{m}_{r-1}, & \text{if } l+1 \le r \le k, \\ 1, & \text{if } r = l. \end{cases}$$

Now, for $l = 0, \ldots, k-1$, we take

$$\widehat{\mathcal{Y}}_k(t_l) = \widehat{\mathcal{M}}_k(t_{l+1}, 1) \cdot \mathcal{G}(t_l, \widehat{X}_k^{\mathrm{WPt}}(t_l))$$

as an approximation to $\mathcal{Y}(t_l)$. Note that, in general, all of the observations $W(t_1), \ldots, W(1)$ are needed to compute the estimate $\widehat{\mathcal{Y}}_k(t_l)$.

Example 2. Consider the linear equation with additive noise from Remark 5. In this case, we have

$$\widehat{\mathcal{M}}_k(t_l, t_r) = \mathcal{M}(t_l, t_r) = 1$$

and

$$\widehat{\mathcal{Y}}_k(t_l) = \mathcal{Y}(t_l) = -\sigma'(t_l).$$

For the linear equation from Example 1, we obtain

$$\widehat{\mathcal{Y}}_k(t_l) = -\beta'(t_l) \cdot \widehat{X}_k^{\mathrm{WPt}}(t_l) \times \prod_{r=l+1}^{k-1} (1 + \alpha(t_r) \cdot (t_{r+1} - t_r) + \beta(t_r) \cdot (W(t_{r+1}) - W(t_r))).$$

4.3. *The basic adaptive scheme* $\widehat{X}_k^\mu$. Choose measurable functions

$$f_l : \mathbb{R}^k \to \mathbb{N}_0$$

for $l = 0, \dots, k-1$. The numbers

$$\mu_l = f_l(\widehat{\mathcal{Y}}_k(t_0), \dots, \widehat{\mathcal{Y}}_k(t_{k-1}))$$

determine the adaptive equidistant discretizations

$$\tau_{l,r} = t_l + r/(k \cdot (\mu_l + 1)), \qquad r = 0, \dots, \mu_l + 1,$$

of the subintervals $[t_l, t_{l+1}]$.

Next, the totality of observations $W(\tau_{l,r})$ is used to estimate the difference $\widehat{X}_k^{\mathrm{WP}} - \widehat{X}_k^{\mathrm{WPt}}$. Put $\mu = (\mu_0, \dots, \mu_{k-1})$, and let $\widetilde{W}^\mu$ denote the piecewise linear interpolation of $W$ at the sites $\tau_{l,r}$. Define the scheme $\widehat{Q}_k^\mu$ by $\widehat{Q}_k^\mu(0) = 0$ and

$$\begin{aligned} \widehat{Q}_k^\mu(t_{l+1}) = \Big(1 &+ a^{(0,1)}(t_l, \widehat{X}_k^{\mathrm{WPt}}(t_l)) \cdot (t_{l+1} - t_l) \\ &+ \sigma^{(0,1)}(t_l, \widehat{X}_k^{\mathrm{WPt}}(t_l)) \cdot (W(t_{l+1}) - W(t_l))\Big) \cdot \widehat{Q}_k^\mu(t_l) \\ &+ \mathcal{G}(t_l, \widehat{X}_k^{\mathrm{WPt}}(t_l)) \cdot \int_{t_l}^{t_{l+1}} (\widetilde{W}^\mu(s) - W(t_l))\, ds \end{aligned}$$

for $l = 0, \dots, k-1$. Note that

$$\widehat{Q}_k^\mu(t_l) = \sum_{r=0}^{l-1} \widehat{\mathcal{Y}}_k(t_r) \cdot \int_{t_r}^{t_{r+1}} (\widetilde{W}^\mu(t) - W(t_r))\, dt. \tag{11}$$

Finally, the basic scheme $\widehat{X}_k^\mu$ is defined by

$$\widehat{X}_k^\mu(t_l) = \widehat{X}_k^{\mathrm{WPt}}(t_l) + \widehat{Q}_k^\mu(t_l), \qquad l = 0, \dots, k.$$

The resulting approximation $\widehat{X}_k^\mu(1)$ belongs to the class $\mathcal{X}^{**}$ and is determined up to the parameters $k$ and $\mu$. Clearly, the number $k$ of the nonadaptive evaluation points should be small compared to the total number $\sum_l \mu_l$ of the adaptively chosen points in order to keep track of the random weight $\mathcal{Y}$. On the other hand, $k$ must be large enough to obtain a sufficiently good approximation $\widehat{\mathcal{Y}}_k$ to $\mathcal{Y}$. Finally, the number $\mu_l$ of observations of $W$ in the interval $]t_l, t_{l+1}[$ should be chosen according to the respective local size of $\mathcal{Y}$. We present three versions of $\widehat{X}_k^\mu(1)$ that are based on this principle.

4.3.1. *The scheme* $\widehat{X}_{p,n}^{**}$ *with varying number of observations of* $W$. Choose a sequence $k_n \in \mathbb{N}$ such that

$$\lim_{n\to\infty} k_n/n = \lim_{n\to\infty} n/k_n^{3/2} = 0 \tag{12}$$

and put

$$\overline{\mathcal{Y}}_{k_n} = \left( \frac{1}{k_n} \sum_{l=0}^{k_n-1} |\widehat{\mathcal{Y}}_{k_n}(t_l)|^{2/3} \right)^{3/2}.$$

Let

$$\mu_l^{(n)} = \begin{cases} \left\lfloor n \cdot \left( |\widehat{\mathcal{Y}}_{k_n}(t_l)|^{2/3} \Big/ \sum_{r=0}^{k_n-1} |\widehat{\mathcal{Y}}_{k_n}(t_r)|^{2/3} \right) \cdot (\overline{\mathcal{Y}}_{k_n})^{p/(p+1)} \right\rfloor, & \text{if } \overline{\mathcal{Y}}_{k_n} > 0, \\ 0, & \text{otherwise,} \end{cases}$$

and define

$$\widehat{X}^{**}_{p,n} = \widehat{X}^{\mu^{(n)}}_{k_n},$$

where

$$\mu^{(n)} = (\mu_0^{(n)}, \ldots, \mu_{k_n-1}^{(n)}).$$

Note that the numbers $\mu_l^{(n)}$ crucially depend on the error parameter $p$. If $p = 2$, then

$$\mu_l^{(n)} = \lfloor n/k_n \cdot |\widehat{\mathcal{Y}}_{k_n}(t_l)|^{2/3} \rfloor.$$

If $p \neq 2$, then all of the approximations $\widehat{\mathcal{Y}}_{k_n}(t_l)$ have to be computed beforehand in order to determine the adaptive discretization.

Finally, we mention that the total number of evaluations of $W$ that are used to obtain the approximation $\widehat{X}^{**}_{p,n}(1)$ is roughly given by $n \cdot S^{p/(p+1)}$, where

$$S = \left( \int_0^1 |\mathcal{Y}(t)|^{2/3}\, dt \right)^{3/2}$$

is the pathwise 2/3-seminorm of the weighting process $\mathcal{Y}$. In general, this quantity depends on the trajectory of $\mathcal{Y}$ so that there is no a priori bound on the computation time available for the user. If all approximations have to be computed in the same amount of time, the following version $\widehat{X}^*_n$ of the basic adaptive scheme can be used. However, note that a price has to be paid for this property; see Theorem 2.

4.3.2. *The scheme $\widehat{X}^*_n$ with fixed number of observations of $W$.* In contrast to the scheme $\widehat{X}^{**}_{p,n}$, the adaptive discretization used by the scheme $\widehat{X}^*_n$ does not depend on the error parameter $p$. Let

$$\mu_l^{(n)} = \begin{cases} \left\lfloor (n - k_n) \cdot |\widehat{\mathcal{Y}}_{k_n}(t_l)|^{2/3} \Big/ \sum_{r=0}^{k_n-1} |\widehat{\mathcal{Y}}_{k_n}(t_r)|^{2/3} \right\rfloor, & \text{if } \overline{\mathcal{Y}}_{k_n} > 0, \\ \lfloor (n - k_n)/k_n \rfloor, & \text{otherwise,} \end{cases}$$

and define $\widehat{X}_n^* = \widehat{X}_{k_n}^{\mu^{(n)}}$, where $\mu^{(n)}$ is determined by $\mu_0^{(n)}, \ldots, \mu_{k_n-1}^{(n)}$.

By definition,

$$n - k_n \le k_n + \sum_{l=0}^{k_n-1} \mu_l^{(n)} \le n$$

holds for the total number of observations, so that the resulting approximation $\widehat{X}_n^*(1)$ belongs to the class $\mathcal{X}_n^*$.

4.3.3. *The scheme $\widehat{X}_n$ with prefixed discretization.* Replacing the quantities $|\widehat{\mathcal{Y}}_{k_n}(t_l)|$ by

$$(E|\mathcal{Y}(t_l)|^2)^{1/2}$$

in the definition of the numbers $\mu_l^{(n)}$ in Section 4.3.2, we obtain the scheme $\widehat{X}_n$, which uses the same discretization for every trajectory of the weighting process $\mathcal{Y}$. The resulting approximation $\widehat{X}_n(1)$ thus belongs to the class $\mathcal{X}_n$. Note that this method requires the computation of the second moments of $\mathcal{Y}$, which might be a difficult task in general.

4.4. *Error analysis of the adaptive schemes.* Now we investigate the asymptotic performance of the approximations $\widehat{X}_{p,n}^{**}(1)$, $\widehat{X}_n^*(1)$ and $\widehat{X}_n(1)$. Additionally, we consider the scheme

$$\widehat{X}_n^{\text{equi}} = \widehat{X}_n^0,$$

which only uses the observations $W(1/n), W(2/n), \ldots, W(1)$ of the driving Brownian motion $W$. Thus, $\widehat{X}_n^{\text{equi}}(1) \in \mathcal{X}_n^{\text{equi}}$. Note that $\widehat{X}_n^{\text{equi}}$ is given by

$$\widehat{X}_n^{\text{equi}}(l/n) = \widehat{X}_n^{\text{WPt}}(l/n) + \frac{1}{2n} \sum_{r=0}^{l-1} \widehat{\mathcal{Y}}_n(t_r) \cdot (W((r+1)/n) - W(r/n))$$

for $l = 0, \ldots, n$.

Recall the constants $C_p^{**}$, $C_p^*$, $C_2$ and $C_p^{\text{equi}}$ from Section 3.2.

THEOREM 2. *The adaptive schemes $\widehat{X}_{p,n}^{**}$, $\widehat{X}_n^*$, $\widehat{X}_n$ and the equidistant scheme $\widehat{X}_n^{\text{equi}}$ satisfy:*

(i) $\lim_{n\to\infty} c(\widehat{X}_{p,n}^{**}(1)) \cdot e_p(\widehat{X}_{p,n}^{**}(1)) = 1/\sqrt{12} \cdot C_p^{**}$,
(ii) $\lim_{n\to\infty} n \cdot e_p(\widehat{X}_n^*(1)) = 1/\sqrt{12} \cdot C_p^*$,
(iii) $\lim_{n\to\infty} n \cdot e_2(\widehat{X}_n(1)) = 1/\sqrt{12} \cdot C_2$,
(iv) $\lim_{n\to\infty} n \cdot e_p(\widehat{X}_n^{\text{equi}}(1)) = 1/\sqrt{12} \cdot C_p^{\text{equi}}$.

Combining Theorem 2 with Theorem 1 from Section 3.2, we immediately obtain

THEOREM 3. *Assume $C_2^{\text{equi}} > 0$. Then the schemes $\widehat{X}_{p,n}^{**}$, $\widehat{X}_n^*$ and $\widehat{X}_n^{\text{equi}}$ are asymptotically optimal for pathwise approximation at $t = 1$ in the respective classes of methods $\mathcal{X}^{**}$, $\mathcal{X}^*$ and $\mathcal{X}^{\text{equi}}$. Moreover, if $p = 2$, then $\widehat{X}_n$ is asymptotically optimal for pathwise approximation at $t = 1$ in the class $\mathcal{X}$.*

REMARK 7. We stress that, in general, the asymptotic constants $C_2/\sqrt{12}$ and $C_2^{\text{equi}}/\sqrt{12}$ cannot be achieved by the Milstein scheme. As an example, consider the equation

$$dX(t) = \sigma(t)\, dW(t), \qquad X(0) = 0,$$

with $\sigma \in C^1([0,1])$. For a discretization

$$0 = t_0 < \cdots < t_n = 1,$$

the corresponding Milstein scheme is given by $\widehat{X}^M_{t_1,\ldots,t_n}(0) = 0$ and

$$\widehat{X}^M_{t_1,\ldots,t_n}(t_l) = \sum_{r=0}^{l-1} \sigma(t_r) \cdot (W(t_{r+1}) - W(t_r)), \qquad l = 1, \ldots, n.$$

Straightforward calculations yield

$$(e_2(\widehat{X}^M_{t_1,\ldots,t_n}(1)))^2 = \tfrac{1}{3} \sum_{l=0}^{n-1} (\sigma'(\xi_l))^2 \cdot (t_{l+1} - t_l)^3 \tag{13}$$

with $t_l \le \xi_l \le t_{l+1}$, so that

$$n^2 \cdot (e_2(\widehat{X}^M_{t_1,\ldots,t_n}(1)))^2 \ge \tfrac{1}{3} \left( \sum_{l=0}^{n-1} (\sigma'(\xi_l))^{2/3} \cdot (t_{l+1} - t_l) \right)^3$$

by the Hölder inequality. Consequently,

$$\liminf_{n\to\infty} n \cdot \inf_{0<t_1<\cdots<t_n=1} e_2(\widehat{X}^M_{t_1,\ldots,t_n}(1)) \ge \frac{1}{\sqrt{3}} \cdot \left( \int_0^1 |\sigma'(t)|^{2/3}\, dt \right)^{3/2} = \frac{C_2}{\sqrt{3}}.$$

Thus, whatever the discretization, the resulting Milstein scheme asymptotically performs suboptimally with respect to pathwise approximation at $t = 1$.

Similarly, for the equidistant Milstein scheme $\widehat{X}_n^M$ we obtain from (13) that

$$\lim_{n\to\infty} n \cdot e_2(\widehat{X}_n^M(1)) = \frac{1}{\sqrt{3}} \cdot \left( \int_0^1 |\sigma'(t)|^2\, dt \right)^{1/2} = \frac{C_2^{\text{equi}}}{\sqrt{3}}.$$

**5. Proofs.** We introduce an auxiliary scheme $\overline{X}_k^{\mathrm{aux}}$ corresponding to the equidistant discretization $t_l = l/k$, $l = 0, \dots, k$, and separately analyze $X(1) - \overline{X}_k^{\mathrm{aux}}(1)$ and $\overline{X}_k^{\mathrm{aux}}(1) - \widehat{X}(1)$ for a method $\widehat{X}(1) \in \mathcal{X}^{**}$. The scheme $\overline{X}_k^{\mathrm{aux}}$ is defined by

$$\overline{X}_k^{\mathrm{aux}}(t_l) = \widehat{X}_k^{\mathrm{WPt}}(t_l) + \overline{Q}_k(t_l), \qquad l = 0, \dots, k.$$

Here $\widehat{X}_k^{\mathrm{WPt}}$ is the truncated Wagner–Platen scheme (see Section 4.1) and the scheme $\overline{Q}_k$ is given by $\overline{Q}_k(0) = 0$ and

$$\begin{aligned}\overline{Q}_k(t_{l+1}) = \Big(1 &+ a^{(0,1)}(t_l, \widehat{X}_k^{\mathrm{WPt}}(t_l)) \cdot (t_{l+1} - t_l) \\ &+ \sigma^{(0,1)}(t_l, \widehat{X}_k^{\mathrm{WPt}}(t_l)) \cdot (W(t_{l+1}) - W(t_l))\Big) \cdot \overline{Q}_k(t_l) \\ &+ \mathcal{G}(t_l, \widehat{X}_k^{\mathrm{WPt}}(t_l)) \cdot \int_{t_l}^{t_{l+1}} (W(s) - W(t_l))\, ds\end{aligned}$$

for $l = 0, \dots, k-1$. Note that

$$\overline{Q}_k(t_l) = \sum_{r=0}^{l-1} \widehat{\mathcal{Y}}_k(t_r) \cdot \int_{t_r}^{t_{r+1}} (W(t) - W(t_l))\, dt. \tag{14}$$

Due to Lemma 12 in the Appendix, we have

$$E|X(1) - \overline{X}_k^{\mathrm{aux}}(1)|^p = O(k^{-3p/2}). \tag{15}$$

Thus, asymptotically $E|\overline{X}_k^{\mathrm{aux}}(1) - \widehat{X}(1)|^p$ will be the dominating term if $k$ is chosen suitable as a function of $c(\widehat{X}(1))$.

We briefly outline the structure of this section. Basic facts on moments of integrated Brownian bridges are stated in Section 5.1. Section 5.2 contains error bounds for the discrete approximation $\widehat{\mathcal{Y}}_k$ of the random weight $\mathcal{Y}$. The lower bounds in Theorem 1 are proven in Section 5.3. The matching upper bounds in Theorem 2 are proven in Section 5.4.

Throughout the following we use $c$ to denote unspecified positive constants that only depend on the error parameter $p$ and the constants from conditions (A) and (B) in Section 2.

5.1. *Moments of integrated Brownian bridges.* Let $B$ denote a Brownian bridge on an interval $[S,T] \subset [0,1]$. Straightforward calculations yield

$$E\left(\int_S^T B(t)\, dt\right)^2 = 1/12 \cdot (T-S)^3. \tag{16}$$

Furthermore, if

$$S = \tau_0 < \cdots < \tau_n = T$$

and $B_1,\dots,B_n$ are independent Brownian bridges on the intervals $[\tau_0,\tau_1],\dots,[\tau_{n-1},\tau_n]$, respectively, then

$$E\left(\sum_{r=0}^{n-1}\int_{\tau_r}^{\tau_{r+1}} B_r(t)\,dt\right)^2 \ge 1/12\cdot(T-S)^3\cdot 1/n^2 \tag{17}$$

by the Hölder inequality.

5.2. *Error bounds for the estimates $\widehat{\mathcal{Y}}_k$.* Recall the discrete approximation $\widehat{\mathcal{M}}_k$ of the field $\mathcal{M}$ from Section 4.2.

LEMMA 1. *For $0\le l\le k-1$, it holds*

$$E|\mathcal{M}(t_l,1)-\widehat{\mathcal{M}}_k(t_l,1)|^{2p}\le c/k^p.$$

PROOF. Note that, by boundedness of $a^{(0,1)}$ and $\sigma^{(0,1)}$,

$$E|\mathcal{M}(t_1,s_1)-\mathcal{M}(t_2,s_2)|^q\le c\cdot c(q)\cdot(\max(|t_1-t_2|,|s_1-s_2|))^{q/2} \tag{18}$$

for all $q\ge 1$, $0\le t_1\le s_1\le 1$, $0\le t_2\le s_2\le 1$, where $c(q)$ only depends on $q$.

Fix $l$ and define the process $\overline{\mathcal{M}}_k(t_l,\cdot)$ on $[t_l,1]$ by $\overline{\mathcal{M}}_k(t_l,t_l)=1$ and

$$\overline{\mathcal{M}}_k(t_l,t)=\overline{\mathcal{M}}_k(t_l,t_r)\cdot\Big(1+a^{(0,1)}(t_r,\widehat{X}_k^{\mathrm{WPt}}(t_r))\cdot(t-t_r) + \sigma^{(0,1)}(t_r,\widehat{X}_k^{\mathrm{WPt}}(t_r))\cdot(W(t)-W(t_r))\Big)$$

for $t\in[t_r,t_{r+1}]$, $r=l,\dots,k-1$. Clearly, $\overline{\mathcal{M}}_k(t_l,t_r)=\widehat{\mathcal{M}}_k(t_l,t_r)$ for $r=l,\dots,k$, and boundedness of $a^{(0,1)}$ and $\sigma^{(0,1)}$ implies

$$E\Big(\sup_{t_l\le t\le 1}|\overline{\mathcal{M}}_k(t_l,t)|^q\Big)\le c\cdot c(q) \tag{19}$$

for every $q\ge 1$, where the constant $c(q)$ only depends on $q$.

Let $t\in[t_l,1]$. Due to (10) we have

$$\begin{aligned}
&|\mathcal{M}(t_l,t)-\overline{\mathcal{M}}_k(t_l,t)|^{2p}\\
&\quad\le c\int_{t_l}^{t}\sum_{r=l}^{k-1}|a^{(0,1)}(s,X(s))\mathcal{M}(t_l,s)\\
&\qquad\qquad - a^{(0,1)}(t_r,\widehat{X}_k^{\mathrm{WPt}}(t_r))\overline{\mathcal{M}}_k(t_l,t_r)|^{2p}\mathbb{1}_{]t_r,t_{r+1}]}(s)\,ds\\
&\quad\quad + c\left|\int_{t_l}^{t}\sum_{r=l}^{k-1}(\sigma^{(0,1)}(s,X(s))\mathcal{M}(t_l,s)\right.\\
&\qquad\qquad \left.-\sigma^{(0,1)}(t_r,\widehat{X}_k^{\mathrm{WPt}}(t_r))\overline{\mathcal{M}}_k(t_l,t_r))\mathbb{1}_{]t_r,t_{r+1}]}(s)\,dW(s)\right|^{2p}.
\end{aligned}$$

Put $V(t)=\sup_{t_l\le s\le t}|\mathcal{M}(t_l,s)-\overline{\mathcal{M}}_k(t_l,s)|$. By the Burkholder inequality,

$$
\begin{aligned}
E|V(t)|^{2p} \le c \int_{t_l}^{t} \sum_{r=l}^{k-1} & E|a^{(0,1)}(s,X(s))\mathcal{M}(t_l,s) \\
& - a^{(0,1)}(t_r,\widehat{X}_k^{\mathrm{WPt}}(t_r))\overline{\mathcal{M}}_k(t_l,t_r)|^{2p} \mathbb{1}_{]t_r,t_{r+1}]}(s)\,ds \\
+ c \int_{t_l}^{t} \sum_{r=l}^{k-1} & E|\sigma^{(0,1)}(s,X(s))\mathcal{M}(t_l,s) \\
& - \sigma^{(0,1)}(t_r,\widehat{X}_k^{\mathrm{WPt}}(t_r))\overline{\mathcal{M}}_k(t_l,t_r)|^{2p} \mathbb{1}_{]t_r,t_{r+1}]}(s)\,ds.
\end{aligned}
$$

Let $s \in [t_r, t_{r+1}]$. By (A),

$$
\begin{aligned}
&|a^{(0,1)}(s,X(s)) \cdot \mathcal{M}(t_l,s) - a^{(0,1)}(t_r,\widehat{X}_k^{\mathrm{WPt}}(t_r)) \cdot \overline{\mathcal{M}}_k(t_l,t_r)| \\
&\quad \le |a^{(0,1)}(s,X(s)) \cdot (\mathcal{M}(t_l,s) - \mathcal{M}(t_l,t_r))| \\
&\qquad + |(a^{(0,1)}(s,X(s)) - a^{(0,1)}(t_r,X(s))) \cdot \mathcal{M}(t_l,t_r)| \\
&\qquad + |(a^{(0,1)}(t_r,X(s)) - a^{(0,1)}(t_r,X(t_r))) \cdot \mathcal{M}(t_l,t_r)| \\
&\qquad + |(a^{(0,1)}(t_r,X(t_r)) - a^{(0,1)}(t_r,\widehat{X}_k^{\mathrm{WPt}}(t_r))) \cdot \mathcal{M}(t_l,t_r)| \\
&\qquad + |a^{(0,1)}(t_r,\widehat{X}_k^{\mathrm{WPt}}(t_r)) \cdot (\mathcal{M}(t_l,t_r) - \overline{\mathcal{M}}_k(t_l,t_r))| \\
&\quad \le c \cdot |\mathcal{M}(t_l,s) - \mathcal{M}(t_l,t_r)| \\
&\qquad + c \cdot ((1 + |X(s)|) \cdot (s - t_r) \\
&\qquad\qquad + |X(s) - X(t_r)| + |X(t_r) - \widehat{X}_k^{\mathrm{WPt}}(t_r)|) \cdot |\mathcal{M}(t_l,t_r)| \\
&\qquad + c \cdot V(s).
\end{aligned}
$$

Observing (5), (6), Lemma 10 and (18), we thus obtain

$$
\begin{aligned}
&E|a^{(0,1)}(s,X(s)) \cdot \mathcal{M}(t_l,s) - a^{(0,1)}(t_r,\widehat{X}_k^{\mathrm{WPt}}(t_r)) \cdot \overline{\mathcal{M}}_k(t_l,t_r)|^{2p} \\
&\quad \le c \cdot (E|X(t_r)|^{4p} + E|X(s) - X(t_r)|^{4p} \\
&\qquad\qquad + E|X(t_r) - \widehat{X}_k^{\mathrm{WPt}}(t_r)|^{4p})^{1/2} \cdot (E|\mathcal{M}(t_l,t_r)|^{4p})^{1/2} \\
&\qquad + c \cdot E|V(s)|^{2p} \\
&\quad \le c/k^p + c \cdot E|V(s)|^{2p},
\end{aligned}
$$

and the same inequality holds with $\sigma^{(0,1)}$ in place of $a^{(0,1)}$.

Consequently, for every $t \in [t_l, 1]$,

$$
E|V(t)|^{2p} \le c/k^p + c \cdot \int_{t_l}^{t} E|V(s)|^{2p}\,ds.
$$

Moreover, by (18) and (19),

$$
E\Big(\sup_{t_l \le t \le 1} |V(t)|^{2p}\Big) < \infty.
$$

Thus, Gronwall's lemma yields

$$\sup_{t_l \le t \le 1} E|V(t)|^{2p} \le c/k^p,$$

which completes the proof. $\square$

LEMMA 2. *For* $0 \le l \le k-1$, *it holds*

$$E|\widehat{\mathcal{Y}}_k(t_l) - \mathcal{Y}(t_l)|^p \le c/k^{p/2}.$$

PROOF. Due to (A),

$$|\mathcal{G}(t_l, \widehat{X}_k^{\mathrm{WPt}}(t_l)) - \mathcal{G}(t_l, X(t_l))| \le (1 + |X(t_l)|^2 + |\widehat{X}_k^{\mathrm{WPt}}(t_l)|^2) \cdot |\widehat{X}_k^{\mathrm{WPt}}(t_l) - X(t_l)|$$

and

$$|\mathcal{G}(t_l, \widehat{X}_k^{\mathrm{WPt}}(t_l))| \le c \cdot (1 + |\widehat{X}_k^{\mathrm{WP}}(t_l)|^2).$$

Hence, by (5), (18), Lemmas 1 and 10,

$$\begin{aligned}
&E|\widehat{\mathcal{Y}}_k(t_l) - \mathcal{Y}(t_l)|^p \\
&\quad \le c \cdot E|(\widehat{\mathcal{M}}_k(t_{l+1}, 1) - \mathcal{M}(t_l, 1)) \cdot \mathcal{G}(t_l, \widehat{X}_k^{\mathrm{WPt}}(t_l))|^p \\
&\qquad + c \cdot E|(\mathcal{M}(t_l, 1) \cdot (\mathcal{G}(t_l, \widehat{X}_k^{\mathrm{WPt}}(t_l)) - \mathcal{G}(t_l, X(t_l))))|^p \\
&\quad \le c \cdot (E|\widehat{\mathcal{M}}_k(t_{l+1}, 1) - \mathcal{M}(t_l, 1)|^{2p})^{1/2} \cdot (E(1 + |\widehat{X}_k^{\mathrm{WPt}}(t_l)|^{4p}))^{1/2} \\
&\qquad + c \cdot (E|\mathcal{M}(t_l, 1)|^{4p})^{1/4} \cdot (E|X(t_l)|^{4p} + E|\widehat{X}_k^{\mathrm{WPt}}(t_l)|^{4p})^{1/2} \\
&\qquad \times (E|\widehat{X}_k^{\mathrm{WPt}}(t_l) - X(t_l)|^{4p})^{1/4} \\
&\quad \le c/k^{p/2},
\end{aligned}$$

which completes the proof. $\square$

Put

$$R = \frac{1}{k} \sum_{l=0}^{k-1} |\mathcal{Y}(t_l)|^{2/3}, \qquad \widehat{R} = \frac{1}{k} \sum_{l=0}^{k-1} |\widehat{\mathcal{Y}}_k(t_l)|^{2/3}.$$

LEMMA 3. *If* $1 \le p \le 2$, *then*

$$|E(\widehat{R}^{3p/2(p+1)}) - E(R^{3p/2(p+1)})| \le c/k^{p/2(p+1)}.$$

*If* $p \ge 2$, *then*

$$|(E(\widehat{R}^{3p/2(p+1)}))^{2(p+1)/3p} - (E(R^{3p/2(p+1)}))^{2(p+1)/3p}| \le c/k^{1/3}.$$

*Furthermore,*

$$|(E(\widehat{R}^{3p/2}))^{2/3p} - (E(R^{3p/2}))^{2/3p}| \le c/k^{1/3}.$$

PROOF. Clearly,

$$|\widehat{R}-R| \le \frac{1}{k}\sum_{l=0}^{k-1} |\widehat{\mathcal{Y}}_k(t_l)-\mathcal{Y}(t_l)|^{2/3}.$$

Assume $1 \le p \le 2$. Then $3p/2(p+1) \le 1$ and we obtain

$$\begin{aligned} E|\widehat{R}^{3p/2(p+1)} - R^{3p/2(p+1)}| &\le E|\widehat{R}-R|^{3p/2(p+1)} \\ &\le E\left(\frac{1}{k}\sum_{l=0}^{k-1} |\widehat{\mathcal{Y}}_k(t_l)-\mathcal{Y}(t_l)|^{2/3}\right)^{3p/2(p+1)} \\ &\le \left(\frac{1}{k}\sum_{l=0}^{k-1} E|\widehat{\mathcal{Y}}_k(t_l)-\mathcal{Y}(t_l)|^{2/3}\right)^{3p/2(p+1)} \\ &\le c/k^{p/2(p+1)} \end{aligned}$$

by Lemma 2, which proves the first inequality.

Next, let $p \ge 2$. Then $3p/2(p+1) \ge 1$. By Lemma 2,

$$\begin{aligned} &|(E(\widehat{R}^{3p/2(p+1)}))^{2(p+1)/3p} - (E(R^{3p/2(p+1)}))^{2(p+1)/3p}| \\ &\quad \le (E|\widehat{R}-R|^{3p/2(p+1)})^{2(p+1)/3p} \\ &\quad \le \frac{1}{k}\sum_{l=0}^{k-1} (E|\widehat{\mathcal{Y}}_k(t_l)-\mathcal{Y}(t_l)|^{p/(p+1)})^{2(p+1)/3p} \\ &\quad \le c/k^{1/3}, \end{aligned}$$

which shows the second inequality. The third inequality is established in the same way. □

5.3. *Proof of the lower bounds in Theorem 1.* Consider an arbitrary sequence of methods $\widehat{X}_N(1) \in \mathcal{X}_N^{**}$. Take a sequence of positive integers $k_N$ that satisfies

$$\lim_{N\to\infty} N/k_N^{3/2} = \lim_{N\to\infty} k_N/N = 0 \tag{20}$$

and assume without loss of generality that $\widehat{X}_N(1)$ uses in particular the evaluation sites

$$t_l = l/k_N, \qquad l = 0,\dots,k_N.$$

Let $d_l^{(N)}$ denote the number of evaluation points that are used by $\widehat{X}_N(1)$ in the interval $]t_l, t_{l+1}[$ and put

$$A_N = \left(\sum_{l=0}^{k_N-1} (\mathcal{Y}(t_l)/(d_l^{(N)}+1))^2\right)^{p/2}.$$

LEMMA 4.

$$\liminf_{N\to\infty} N \cdot e_p(\widehat{X}_N(1)) \geq \mathfrak{m}_p/12^{1/2} \cdot \liminf_{N\to\infty} N/k_N^{3/2} \cdot (E(A_N))^{1/p}.$$

PROOF. By (15),

$$e_p(\widehat{X}_N(1)) \geq (E|\overline{X}^{\mathrm{aux}}_{k_N}(1) - \widehat{X}_N(1)|^p)^{1/p} - c/k_N^{3/2}. \tag{21}$$

Let $\mathfrak{A}_N$ denote the $\sigma$-algebra that is generated by the data used by $\widehat{X}_N(1)$ and put $Z = W - E(W|\mathfrak{A}_N)$ as well as

$$V = \widehat{X}_N(1) - \widehat{X}^{\mathrm{WPt}}_{k_N}(1) - \sum_{l=0}^{k_N-1} \widehat{\mathcal{Y}}_{k_N}(t_l) \cdot \int_{t_l}^{t_{l+1}} (E(W(t)|\mathfrak{A}_N) - W(t_l))\, dt.$$

By definition of $\overline{X}^{\mathrm{aux}}_{k_N}$ and (14), we have

$$\overline{X}^{\mathrm{aux}}_{k_N}(1) - \widehat{X}_N(1) = \sum_{l=0}^{k_N-1} \widehat{\mathcal{Y}}_{k_N}(t_l) \cdot \int_{t_l}^{t_{l+1}} Z(t)\, dt - V.$$

Note that $V$ and the numbers $d_l^{(N)}$ are $\mathfrak{A}_N$-measurable. Conditioned on $\mathfrak{A}_N$, the evaluation sites used by $\widehat{X}_N(1)$ are fixed and the process $Z$ consists of independent Brownian bridges corresponding to the respective subintervals. Hence, by (17),

$$\begin{aligned}
&E(|\overline{X}^{\mathrm{aux}}_{k_N}(1) - \widehat{X}_N(1)|^p|\mathfrak{A}_N) \\
&\quad \geq E\left(\left|\sum_{l=0}^{k_N-1} \widehat{\mathcal{Y}}_{k_N}(t_l) \cdot \int_{t_l}^{t_{l+1}} Z(t)\, dt\right|^p \middle| \mathfrak{A}_N\right) \\
&\quad = \mathfrak{m}_p^p \cdot \left(E\left(\left|\sum_{l=0}^{k_N-1} \widehat{\mathcal{Y}}_{k_N}(t_l) \cdot \int_{t_l}^{t_{l+1}} Z(t)\, dt\right|^2 \middle| \mathfrak{A}_N\right)\right)^{p/2} \qquad (22) \\
&\quad = \mathfrak{m}_p^p \cdot \left(\sum_{l=0}^{k_N-1} (\widehat{\mathcal{Y}}_{k_N}(t_l))^2 \cdot E\left(\left(\int_{t_l}^{t_{l+1}} Z(t)\, dt\right)^2 \middle| \mathfrak{A}_N\right)\right)^{p/2} \\
&\quad \geq \mathfrak{m}_p^p \cdot \left(\sum_{l=0}^{k_N-1} (\widehat{\mathcal{Y}}_{k_N}(t_l))^2 \cdot (12 k_N^3 \cdot (d_l^{(N)} + 1)^2)^{-1}\right)^{p/2} \\
&\quad = \mathfrak{m}_p^p / 12^{p/2} \cdot 1/k_N^{3p/2} \cdot \left(\sum_{l=0}^{k_N-1} (\widehat{\mathcal{Y}}_{k_N}(t_l)/(d_l^{(N)} + 1))^2\right)^{p/2}.
\end{aligned}$$

Combine (21) with (22) to obtain

$$\liminf_{N\to\infty} N \cdot e_p(\widehat{X}_N(1))$$

$$\geq \frac{\mathfrak{m}_p}{\sqrt{12}} \cdot \liminf_{N\to\infty} N/k_N^{3/2} \cdot \left( E\left( \sum_{l=0}^{k_N-1} (\widehat{\mathcal{Y}}_{k_N}(t_l)/(d_l^{(N)}+1))^2 \right)^{p/2} \right)^{1/p}.$$

Let $q=\max(2,p)$. Lemma 2 implies

$$\begin{aligned}
&\left| \left( E\left( \sum_{l=0}^{k_N-1} (\widehat{\mathcal{Y}}_{k_N}(t_l)/(d_l^{(N)}+1))^2 \right)^{p/2} \right)^{1/p} - (E(A_N))^{1/p} \right| \\
&\quad \leq \left( E\left( \sum_{l=0}^{k_N-1} |\widehat{\mathcal{Y}}_{k_N}(t_l) - \mathcal{Y}(t_l)|^2 \right)^{p/2} \right)^{1/p} \\
&\quad \leq \left( E\left( \sum_{l=0}^{k_N-1} |\widehat{\mathcal{Y}}_{k_N}(t_l) - \mathcal{Y}(t_l)|^2 \right)^{q/2} \right)^{1/q} \\
&\quad \leq \left( \sum_{l=0}^{k_N-1} (E|\widehat{\mathcal{Y}}_{k_N}(t_l) - \mathcal{Y}(t_l)|^q)^{2/q} \right)^{1/q} \\
&\quad \leq c.
\end{aligned}$$

Thus, by (20),

$$\begin{aligned}
&\liminf_{N\to\infty} N/k_N^{3/2} \cdot \left( E\left( \sum_{l=0}^{k_N-1} (\widehat{\mathcal{Y}}_{k_N}(t_l)/(d_l^{(N)}+1))^2 \right)^{p/2} \right)^{1/p} \\
&\quad \geq \liminf_{N\to\infty} N/k_N^{3/2} \cdot (E(A_N))^{1/p},
\end{aligned}$$

which completes the proof. □

Now, we analyze the classes $\mathcal{X}^{**}$, $\mathcal{X}^{*}$, $\mathcal{X}^{\text{equi}}$ and the class $\mathcal{X}$ in the case $p=2$.

LEMMA 5.

(i) *If $\widehat{X}_N(1) \in \mathcal{X}_N^{**}$ for every $N$, then*

$$\liminf_{N\to\infty} N/k_N^{3/2} \cdot (E(A_N))^{1/p} \geq C_p^{**}/\mathfrak{m}_p.$$

(ii) *If $\widehat{X}_N(1) \in \mathcal{X}_N^{*}$ for every $N$, then*

$$\liminf_{N\to\infty} N/k_N^{3/2} \cdot (E(A_N))^{1/p} \geq C_p^{*}/\mathfrak{m}_p.$$

(iii) *If $p=2$ and $\widehat{X}_N(1) \in \mathcal{X}_N$ for every $N$, then*

$$\liminf_{N\to\infty} N/k_N^{3/2} \cdot (E(A_N))^{1/2} \geq C_2.$$

(iv) *If* $\widehat{X}_N(1) \in \mathcal{X}_N^{\mathrm{equi}}$ *for every* $N$, *then*

$$\liminf_{N\to\infty} \; N/k_N^{3/2} \cdot (E(A_N))^{1/p} \geq C_p^{\mathrm{equi}}/\mathfrak{m}_p.$$

PROOF. By definition of $\mathcal{X}_N^{**}$ and the Hölder inequality,

$$\begin{aligned}
N^{p/(p+1)} \cdot (E(A_N))^{1/(p+1)} &\\
&\geq \left(E\sum_{l=0}^{k_N-1}(d_l^{(N)}+1)\right)^{p/(p+1)} \cdot (E(A_N))^{1/(p+1)}\\
&\geq E\left(\left(\sum_{l=0}^{k_N-1}(d_l^{(N)}+1)\right)^{p/(p+1)} \cdot A_N^{1/(p+1)}\right)\\
&= E\left(\left(\sum_{l=0}^{k_N-1}(d_l^{(N)}+1)\right)^{2/3} \cdot \left(\sum_{l=0}^{k_N-1}(\mathcal{Y}(t_l)/(d_l^{(N)}+1))^2\right)^{1/3}\right)^{3p/(2(p+1))}\\
&\geq E\left(\sum_{l=0}^{k_N-1}|\mathcal{Y}(t_l)|^{2/3}\right)^{3p/(2(p+1))}.
\end{aligned}$$

Hence (i) follows from

$$\begin{aligned}
\liminf_{N\to\infty} N/k_N^{3/2} \cdot (E(A_N))^{1/p} &\\
&\geq \liminf_{N\to\infty}\left(E\left(\frac{1}{k_N}\sum_{l=0}^{k_N-1}|\mathcal{Y}(t_l)|^{2/3}\right)^{3p/(2(p+1))}\right)^{(p+1)/p}\\
&\geq \left(E\left(\liminf_{N\to\infty}\frac{1}{k_N}\sum_{l=0}^{k_N-1}|\mathcal{Y}(t_l)|^{2/3}\right)^{3p/(2(p+1))}\right)^{(p+1)/p}\\
&= C_p^{**}/\mathfrak{m}_p.
\end{aligned}$$

By definition of $\mathcal{X}_N^{*}$,

$$\begin{aligned}
N^p \cdot E(A_N) &\geq E\left(\left(\sum_{l=0}^{k_N-1}(d_l^{(N)}+1)\right)^2 \cdot \left(\sum_{l=0}^{k_N-1}(\mathcal{Y}(t_l)/(d_l^{(N)}+1))^2\right)\right)^{p/2}\\
&\geq E\left(\sum_{l=0}^{k_N-1}|\mathcal{Y}(t_l)|^{2/3}\right)^{3p/2},
\end{aligned}$$

so that

$$\liminf_{N\to\infty} \; N/k_N^{3/2} \cdot (E(A_N))^{1/p}$$

$$
\begin{aligned}
&\geq \liminf_{N\to\infty} \left( E\left( \frac{1}{k_N} \sum_{l=0}^{k_N-1} |\mathcal{Y}(t_l)|^{2/3} \right)^{3p/2} \right)^{1/p} \\
&\geq \left( E\left( \liminf_{N\to\infty} \frac{1}{k_N} \sum_{l=0}^{k_N-1} |\mathcal{Y}(t_l)|^{2/3} \right)^{3p/2} \right)^{1/p} \\
&= C_p^* / \mathfrak{m}_p,
\end{aligned}
$$

which proves (ii).

Next, assume $p=2$. By definition of $\mathcal{X}_N$, the numbers $d_l^{(N)}$ are deterministic. Thus

$$
\begin{aligned}
N^2 \cdot E(A_N) &\geq \left( \sum_{l=0}^{k_N-1} (d_l^{(N)}+1) \right)^2 \cdot \left( \sum_{l=0}^{k_N-1} E(\mathcal{Y}(t_l))^2 / (d_l^{(N)}+1)^2 \right) \\
&\geq \left( \sum_{l=0}^{k_N-1} (E(\mathcal{Y}(t_l))^2)^{1/3} \right)^3 .
\end{aligned}
$$

It follows that

$$
\liminf_{N\to\infty} N/k_N^{3/2} \cdot (E(A_N))^{1/2} \geq \liminf_{N\to\infty} \left( \frac{1}{k_N} \sum_{l=0}^{k_N-1} (E(\mathcal{Y}(t_l))^2)^{1/3} \right)^{3/2} = C_2,
$$

which shows (iii).

Finally, by definition of $\mathcal{X}_N^{\mathrm{equi}}$, the numbers $d_l^{(N)}$ are deterministic with

$$
d_0^{(N)} = \cdots = d_{k_N-1}^{(N)}.
$$

Hence,

$$
N \cdot (E(A_N))^{1/p} \geq k_N \left( E\left( \sum_{l=0}^{k_N-1} (\mathcal{Y}(t_l))^2 \right)^{p/2} \right)^{1/p}.
$$

Consequently,

$$
\begin{aligned}
\liminf_{N\to\infty} N/k_N^{3/2} \cdot (E(A_N))^{1/p} &\geq \liminf_{N\to\infty} \left( E\left( \frac{1}{k_N} \sum_{l=0}^{k_N-1} (\mathcal{Y}(t_l))^2 \right)^{p/2} \right)^{1/p} \\
&\geq \left( E\left( \liminf_{N\to\infty} \frac{1}{k_N} \sum_{l=0}^{k_N-1} (\mathcal{Y}(t_l))^2 \right)^{p/2} \right)^{1/p} \\
&= C_p^{\mathrm{equi}} / \mathfrak{m}_p,
\end{aligned}
$$

which completes the proof. □

Combine Lemma 4 with Lemma 5 to obtain the lower bounds in Theorem 1. Clearly, these lower bounds yield the lower bounds from Theorem 2.

5.4. *Proof of the upper bounds in Theorem 2.* Let $k \in \mathbb{N}$ and consider a basic scheme $\widehat{X}_k^\mu$; see Section 4.3. Put

$$B_k = \left( \sum_{l=0}^{k-1} (\widehat{\mathcal{Y}}_k(t_l)/(\mu_l + 1))^2 \right)^{p/2}.$$

LEMMA 6.

$$e_p(\widehat{X}_k^\mu(1)) \le \mathfrak{m}_p/12^{1/2} \cdot 1/k^{3/2} \cdot (E(B_k))^{1/p} + c/k^{3/2}.$$

PROOF. Due to (15), we have

$$e_p(\widehat{X}_k^\mu(1)) \le (E|\overline{X}_k^{\text{aux}}(1) - \widehat{X}_k^\mu(1)|^p)^{1/p} + c/k^{3/2}. \tag{23}$$

By (11) and (14),

$$\begin{aligned} \overline{X}_k^{\text{aux}}(1) - \widehat{X}_k^\mu(1) &= \overline{Q}_k(1) - \widehat{Q}_k(1) \\ &= \sum_{l=0}^{k-1} \widehat{\mathcal{Y}}_k(t_l) \cdot \int_{t_l}^{t_{l+1}} (W(t) - \widetilde{W}^\mu(t))\, dt. \end{aligned}$$

Let $\mathfrak{B}$ denote the $\sigma$-algebra that is generated by $X(0), W(t_1), \ldots, W(1)$, and recall that the adaptive discretization determined by $\mu$ consists of the $\mathfrak{B}$-measurable points

$$\tau_{l,r} = t_l + r/(k \cdot (\mu_l + 1)), \qquad r = 0, \ldots, \mu_l + 1.$$

Conditioned on $\mathfrak{B}$, the discretization is fixed and the process $W - \widetilde{W}^\mu$ consists of independent Brownian bridges corresponding to the respective subintervals. Using (16), we thus obtain

$$\begin{aligned} &E(|\overline{X}_k^{\text{aux}}(1) - \widehat{X}_k^\mu(1)|^p|\mathfrak{B}) \\ &\quad = \mathfrak{m}_p^p/12^{p/2} \cdot 1/k^{3p/2} \cdot \left( \sum_{l=0}^{k-1} (\widehat{\mathcal{Y}}_k(t_l)/(\mu_l + 1))^2 \right)^{p/2}. \end{aligned} \tag{24}$$

Combine (23) with (24) to obtain the desired result. □

Now we turn to the specific schemes $\widehat{X}_{p,n}^{**}$, $\widehat{X}_n^*$, $\widehat{X}_n$ and $\widehat{X}_n^{\text{equi}}$.

LEMMA 7. *The scheme $\widehat{X}_{p,n}^{**}$ satisfies*

$$\limsup_{n \to \infty} (c(\widehat{X}_{p,n}^{**}(1)))/n \le E\left( \int_0^1 |\mathcal{Y}(t)|^{2/3}\, dt \right)^{3p/2(p+1)}$$

*and*

$$\limsup_{n \to \infty} n \cdot e_p(\widehat{X}_{p,n}^{**}(1)) \le \mathfrak{m}_p/12^{1/2} \cdot \left( E\left( \int_0^1 |\mathcal{Y}(t)|^{2/3}\, dt \right)^{3p/2(p+1)} \right)^{1/p}.$$

PROOF. By definition,

$$c(\widehat{X}^{**}_{p,n}(1)) \le k_n + n \cdot E\left(\frac{1}{k_n}\sum_{l=0}^{k_n-1} |\widehat{\mathcal{Y}}_k(t_l)|^{2/3}\right)^{3p/2(p+1)}.$$

Observe (12) and use Lemma 3 to get

$$\begin{aligned}\limsup_{n\to\infty}(c(\widehat{X}^{**}_{p,n}(1)))/n &\le \limsup_{n\to\infty} E\left(\frac{1}{k_n}\sum_{l=0}^{k_n-1} |\widehat{\mathcal{Y}}_k(t_l)|^{2/3}\right)^{3p/2(p+1)} \\ &\le E\left(\limsup_{n\to\infty}\frac{1}{k_n}\sum_{l=0}^{k_n-1} |\mathcal{Y}(t_l)|^{2/3}\right)^{3p/2(p+1)} \\ &= E\left(\int_0^1 |\mathcal{Y}(t)|^{2/3}\,dt\right)^{3p/2(p+1)},\end{aligned}$$

which proves the first inequality.

Next, observe that, for the scheme $\widehat{X}^{**}_{p,n}$,

$$\sum_{l=0}^{k_n-1} (\widehat{\mathcal{Y}}_k(t_l)/(\mu_l^{(n)}+1))^2 \le n^2 \cdot k_n^{3p/(p+1)} \cdot \left(\sum_{l=0}^{k_n-1} |\widehat{\mathcal{Y}}_k(t_l)|^{2/3}\right)^{3/(p+1)}.$$

Hence,

$$1/k_n^{3p/2} \cdot B_{k_n} \le 1/n^p \cdot \left(\frac{1}{k_n}\sum_{l=0}^{k_n-1} |\widehat{\mathcal{Y}}_k(t_l)|^{2/3}\right)^{3p/2(p+1)}.$$

Using Lemmas 6 and 3, we thus conclude that

$$\begin{aligned}\limsup_{n\to\infty} n \cdot e_p(\widehat{X}^{**}_{p,n}(1)) &\le \mathfrak{m}_p/12^{1/2} \cdot \limsup_{n\to\infty}\left(E\left(\frac{1}{k_n}\sum_{l=0}^{k_n-1} |\mathcal{Y}(t_l)|^{2/3}\right)^{3p/2(p+1)}\right)^{1/p} \\ &\le \mathfrak{m}_p/12^{1/2} \cdot \left(E\left(\int_0^1 |\mathcal{Y}(t)|^{2/3}\,dt\right)^{3p/2(p+1)}\right)^{1/p},\end{aligned}$$

which completes the proof. □

Clearly, Lemma 7 implies the upper bound in Theorem 2(i).

LEMMA 8. *The schemes* $\widehat{X}^*_n$ *and* $\widehat{X}_n$ *satisfy*

$$\limsup_{n\to\infty} n \cdot e_p(\widehat{X}^*_n(1)) \le C^*_p/\sqrt{12}$$

*and*

$$\limsup_{n\to\infty} n \cdot e_2(\widehat{X}_n(1)) \le C_2/\sqrt{12}.$$

PROOF. By definition of $\widehat{X}_n^*(1)$,

$$\sum_{l=0}^{k_n-1} (\widehat{\mathcal{Y}}_k(t_l)/(\mu_l^{(n)}+1))^2 \le 1/(n-k_n)^2 \cdot \left(\sum_{l=0}^{k_n-1} |\widehat{\mathcal{Y}}_k(t_l)|^{2/3}\right)^3.$$

Thus,

$$1/k_n^{3p/2} \cdot B_{k_n} \le 1/(n-k_n)^p \cdot \left(\sum_{l=0}^{k_n-1} |\widehat{\mathcal{Y}}_k(t_l)|^{2/3}\right)^{3p/2}.$$

Use Lemmas 6 and 3 to obtain

$$\begin{aligned}\limsup_{n\to\infty} n \cdot e_p(\widehat{X}_n^*(1)) &\le \mathfrak{m}_p/12^{1/2} \cdot \limsup_{n\to\infty}\left(E\left(\frac{1}{k_n}\sum_{l=0}^{k_n-1} |\mathcal{Y}(t_l)|^{2/3}\right)^{3p/2}\right)^{1/p} \\ &\le \mathfrak{m}_p/12^{1/2} \cdot \left(E\left(\int_0^1 |\mathcal{Y}(t)|^{2/3}\,dt\right)^{3p/2}\right)^{1/p},\end{aligned}$$

which establishes the first inequality.

By definition of $\widehat{X}_n$, the numbers $\mu_l^{(n)}$ are deterministic with

$$\sum_{l=0}^{k_n-1} E|\mathcal{Y}(t_l)|^2/(\mu_l^{(n)}+1)^2 \le 1/(n-k_n)^2 \cdot \left(\sum_{l=0}^{k_n-1} (E|\mathcal{Y}(t_l)|^2)^{1/3}\right)^3.$$

Furthermore, Lemma 2 implies

$$\left(\sum_{l=0}^{k_n-1} \frac{E|\widehat{\mathcal{Y}}_k(t_l)|^2}{(\mu_l^{(n)}+1)^2}\right)^{1/2} \le \left(\sum_{l=0}^{k_n-1} \frac{E|\mathcal{Y}(t_l)|^2}{(\mu_l^{(n)}+1)^2}\right)^{1/2} + c.$$

Hence, by Lemma 6,

$$n \cdot e_2(\widehat{X}_n(1)) \le \mathfrak{m}_p/12^{1/2} \cdot n/(n-k_n) \cdot \left(\frac{1}{k_n}\sum_{l=0}^{k_n-1} (E|\mathcal{Y}(t_l)|^2)^{1/3}\right)^3 + c \cdot n/k_n^{3/2}.$$

Observing (12), we get

$$\begin{aligned}\limsup_{n\to\infty} n \cdot e_2(\widehat{X}_n(1)) &\le \mathfrak{m}_p/12^{1/2} \cdot \limsup_{n\to\infty}\left(\frac{1}{k_n}\sum_{l=0}^{k_n-1} (E|\mathcal{Y}(t_l)|^2)^{1/3}\right)^{3/2} \\ &= \mathfrak{m}_p/12^{1/2} \cdot \left(\int_0^1 (E|\mathcal{Y}(t)|^2)^{1/3}\right)^{3/2},\end{aligned}$$

which completes the proof. □

Lemma 8 yields the upper bounds in Theorem 2(ii), (iii). It remains to establish the upper bound in Theorem 2(iv).

LEMMA 9. *The scheme* $\widehat{X}_n^{\mathrm{equi}}$ *satisfies*

$$\limsup_{n\to\infty} n \cdot e_p(\widehat{X}_n^{\mathrm{equi}}(1)) \le C_p^{\mathrm{equi}}/\sqrt{12}.$$

PROOF. Lemma 6 yields

$$e_p(\widehat{X}_n^{\mathrm{equi}}(1)) \le \mathfrak{m}_p/12^{1/2} \cdot 1/n^{3/2} \cdot \left( E\left( \sum_{l=0}^{n-1} |\widehat{\mathcal{Y}}_n(t_l)|^2 \right)^{p/2} \right)^{1/p} + c/n^{3/2}.$$

Hence, by Lemma 2,

$$n \cdot e_p(\widehat{X}_n^{\mathrm{equi}}(1)) \le \mathfrak{m}_p/12^{1/2} \cdot \left( E\left( \frac{1}{n} \sum_{l=0}^{n-1} |\mathcal{Y}(l/n)|^2 \right)^{p/2} \right)^{1/p} + c/n^{1/2}.$$

We conclude

$$\begin{aligned} &\limsup_{n\to\infty} n \cdot e_p(\widehat{X}_n^{\mathrm{equi}}(1)) \\ &\qquad \le \mathfrak{m}_p/12^{1/2} \cdot \limsup_{n\to\infty} \left( E\left( \frac{1}{n} \sum_{l=0}^{n-1} |\mathcal{Y}(l/n)|^2 \right)^{p/2} \right)^{1/p} \\ &\qquad \le \mathfrak{m}_p/12^{1/2} \cdot \left( E\left( \int_0^1 |\mathcal{Y}(t)|^2 \, dt \right)^{p/2} \right)^{1/p}, \end{aligned}$$

which completes the proof. □

The upper bounds from Theorem 2 imply the upper bounds from Theorem 1.

## APPENDIX

The goal of this appendix is to establish the error bound (15) for the auxiliary scheme $\overline{X}_k^{\mathrm{aux}}$ from Section 5. Throughout, we fix a discretization

$$0 = t_0 < \cdots < t_k = 1,$$

and we put

$$\Delta_l = t_{l+1} - t_l, \qquad \Delta_{\max} = \max_{l=0,\ldots,k-1} \Delta_l.$$

Moreover, we use $\mathcal{F}_t$ to denote the $\sigma$-algebra that is generated by $X(0)$ and $W(s)$, $0 \le s \le t$. Finally, we put

$$|||Y|||_q = (E|Y|^q)^{1/q}$$

for a random variable $Y$ and $q \ge 1$.

We start with error bounds on continuous versions of the Wagner–Platen scheme and its truncated version. Define processes $X^{\mathrm{WPt}}$ and $X^{\mathrm{WP}}$ by $X^{\mathrm{WPt}}(0) = X^{\mathrm{WP}}(0) = X(0)$,

$$\begin{aligned} X^{\mathrm{WPt}}(t) = {} & X^{\mathrm{WPt}}(t_l) + a(t_l, X^{\mathrm{WPt}}(t_l)) \cdot (t - t_l) + \sigma(t_l, X^{\mathrm{WPt}}(t_l)) \cdot (W(t) - W(t_l)) \\ & + 1/2 \cdot (\sigma\sigma^{(0,1)})(t_l, X^{\mathrm{WPt}}(t_l)) \cdot ((W(t) - W(t_l))^2 - (t - t_l)) \\ & + (\sigma^{(1,0)} + a\sigma^{(0,1)} - 1/2 \cdot \sigma(\sigma^{(0,1)})^2) \\ & \quad \times (t_l, X^{\mathrm{WPt}}(t_l)) \cdot (W(t) - W(t_l)) \cdot (t - t_l) \\ & + 1/6 \cdot (\sigma(\sigma^{(0,1)})^2 + \sigma^2\sigma^{(0,2)})(t_l, X^{\mathrm{WPt}}(t_l)) \cdot (W(t) - W(t_l))^3 \\ & + 1/2 \cdot (a^{(1,0)} + aa^{(0,1)} + 1/2 \cdot \sigma^2 a^{(0,2)})(t_l, X^{\mathrm{WPt}}(t_l)) \cdot (t - t_l)^2 \end{aligned}$$

and

$$\begin{aligned} X^{\mathrm{WP}}(t) = {} & X^{\mathrm{WP}}(t_l) + a(t_l, X^{\mathrm{WP}}(t_l)) \cdot (t - t_l) + \sigma(t_l, X^{\mathrm{WP}}(t_l)) \cdot (W(t) - W(t_l)) \\ & + 1/2 \cdot (\sigma\sigma^{(0,1)})(t_l, X^{\mathrm{WP}}(t_l)) \cdot ((W(t) - W(t_l))^2 - (t - t_l)) \\ & + (\sigma^{(1,0)} + a\sigma^{(0,1)} - 1/2 \cdot \sigma(\sigma^{(0,1)})^2) \\ & \quad \times (t_l, X^{\mathrm{WP}}(t_l)) \cdot (W(t) - W(t_l)) \cdot (t - t_l) \\ & + 1/6 \cdot (\sigma(\sigma^{(0,1)})^2 + \sigma^2\sigma^{(0,2)})(t_l, X^{\mathrm{WP}}(t_l)) \cdot (W(t) - W(t_l))^3 \\ & + 1/2 \cdot (a^{(1,0)} + aa^{(0,1)} + 1/2 \cdot \sigma^2 a^{(0,2)})(t_l, X^{\mathrm{WP}}(t_l)) \cdot (t - t_l)^2 \\ & + \mathcal{G}(t_l, X^{\mathrm{WP}}(t_l)) \cdot \int_{t_l}^{t} (W(s) - W(t_l))\, ds \end{aligned}$$

for $t \in [t_l, t_{l+1}]$, $l = 0, \ldots, k-1$.

LEMMA 10. *The processes $X^{\mathrm{WPt}}$ and $X^{\mathrm{WP}}$ satisfy:*

(i) $$\sup_{t\in[0,1]} E|X^{\mathrm{WPt}}(t)|^{16p} \le c,$$

(ii) $$\sup_{t\in[0,1]} E|X^{\mathrm{WP}}(t)|^{16p} \le c,$$

*as well as*

(iii) $$\sup_{t\in[0,1]} E|X(t) - X^{\mathrm{WPt}}(t)|^{4p} \le c \cdot \Delta_{\max}^{4p},$$

(iv) $$\sup_{t\in[0,1]} E|X(t) - X^{\mathrm{WP}}(t)|^{4p} \le c \cdot \Delta_{\max}^{6p}.$$

See Kloeden and Platen (1995) for a proof of (ii) and (iv) under much stronger assumptions on $a$ and $\sigma$ than stated in (A) in Section 2. For a proof of Lemma 10 under condition (A), we refer to Müller-Gronbach (2002b).

Next, define the process $Q$ by $Q(0) = 0$ and

$$Q(t) = \Big(1 + a^{(0,1)}(t_l, X^{\mathrm{WPt}}(t_l)) \cdot (t - t_l) + \sigma^{(0,1)}(t_l, X^{\mathrm{WPt}}(t_l)) \cdot (W(t) - W(t_l))\Big) \cdot Q(t_l) + \mathcal{G}(t_l, X^{\mathrm{WPt}}(t_l)) \cdot \int_{t_l}^{t} (W(s) - W(t_l))\, ds$$

for $t \in [t_l, t_{l+1}]$. Note that $Q(t_l) = \overline{Q}_k(t_l)$ for an equidistant discretization (9).

Lemma 11. *The process $Q$ satisfies*

$$\sup_{t\in[0,1]} E|Q(t)|^{4p} \le c \cdot \Delta_{\max}^{8p}$$

*and*

$$\sup_{t\in[t_l,t_{l+1}]} E|Q(t) - Q(t_l)|^{4p} \le c \cdot \Delta_{\max}^{6p}.$$

Proof. Fix $t \in [t_l, t_{l+1}]$ and let

$$U = (1 + a^{(0,1)}(t_l, X^{\mathrm{WPt}}(t_l)) \cdot (t - t_l)) \cdot Q(t_l), \qquad V = Q(t) - U.$$

Put $q = \lceil 2p \rceil$ and note that $4p \le 2q \le 8p$. Let $r \in \{1, \dots, 2q\}$. Observing (A), we have

$$E(|V|^r | \mathcal{F}_{t_l}) \le c \cdot |Q(t_l)|^r \cdot (t - t_l)^{r/2} + c \cdot (1 + |X^{\mathrm{WPt}}(t_l)|^r) \cdot (t - t_l)^{3r/2} \tag{25}$$

as well as

$$|U|^r \le (1 + c \cdot (t - t_l)) \cdot |Q(t_l)|^r.$$

Moreover, if $r$ is odd, then

$$E(V^r | \mathcal{F}_{t_l}) = 0.$$

Hence,

$$\begin{aligned} E((Q(t))^{2q} | \mathcal{F}_{t_l}) &= U^{2q} + \sum_{r=1}^{2q} \binom{2q}{r} \cdot U^{2q-r} \cdot E(V^r | \mathcal{F}_{t_l}) \\ &\le (1 + c \cdot (t - t_l)) \cdot |Q(t_l)|^{2q} \\ &\quad + c \cdot \sum_{r=1}^{q} \binom{2q}{2r} \cdot |Q(t_l)|^{2q-2r} \cdot (1 + |X^{\mathrm{WPt}}(t_l)|^{2r}) \cdot (t - t_l)^{3r}. \end{aligned}$$

Use Lemma 10(i) to obtain

$$\begin{aligned}\|\!|Q(t)\|\!|_{2q}^{2q} &\le (1+c\cdot(t-t_l))\cdot \|\!|Q(t_l)\|\!|_{2q}^{2q} + c\cdot\sum_{r=1}^{q}\binom{2q}{2r}\cdot\|\!|Q(t_l)\|\!|_{2q}^{2q-2r}\cdot(t-t_l)^{3r}\\ &\le (1+c\cdot(t-t_l))\cdot\|\!|Q(t_l)\|\!|_{2q}^{2q} + c\cdot(t-t_l)\cdot(\|\!|Q(t_l)\|\!|_{2q}+(t-t_l))^{2q}\\ &\le (1+c\cdot(t-t_l))\cdot\|\!|Q(t_l)\|\!|_{2q}^{2q} + c\cdot(t-t_l)^{2q+1},\end{aligned}$$

so that the first inequality follows from Gronwall's lemma.

Due to (25) and Lemma 10(i),

$$E|V|^{2q} \le c\cdot E|Q(t_l)|^{2q}\cdot(t-t_l)^q + c\cdot(t-t_l)^{3q}.$$

Thus, by (A) and the first inequality,

$$E|Q(t)-Q(t_l)|^{2q} \le c\cdot(t-t_l)^{2q}\cdot E|Q(t_l)|^{2q} + c\cdot E|V|^{2q} \le c\cdot\Delta_{\max}^{3q},$$

which proves the second inequality. □

Finally, we consider the process

$$X^{\mathrm{aux}} = X^{\mathrm{WPt}} + Q.$$

Lemma 12. *The process* $X^{\mathrm{aux}}$ *satisfies*

$$\sup_{t\in[0,1]} E|X(t)-X^{\mathrm{aux}}(t)|^p \le c\cdot\Delta_{\max}^{3p/2}.$$

Note that $X^{\mathrm{aux}}(t_l) = \overline{X}_k^{\mathrm{aux}}(t_l)$ for an equidistant discretization (9). Consequently, Lemma 12 immediately implies (15).

Proof of Lemma 12. In view of Lemma 10(iv), it is enough to show

$$\sup_{t\in[0,1]} E|X^{\mathrm{WP}}(t)-X^{\mathrm{aux}}(t)|^p \le c\cdot\Delta_{\max}^{3p/2}. \tag{26}$$

Let

$$\begin{aligned} g_1 &= 1/2\sigma\sigma^{(0,1)},\\ g_2 &= \sigma^{(1,0)} + a\sigma^{(0,1)} - 1/2\sigma(\sigma^{(0,1)})^2,\\ g_3 &= 1/6(\sigma(\sigma^{(0,1)})^2 + \sigma^2\sigma^{(0,2)}),\\ g_4 &= 1/2(a^{(1,0)} + aa^{(0,1)} + 1/2\sigma^2 a^{(0,2)}),\\ g_5 &= \mathcal{G}.\end{aligned}$$

Fix $t \in [t_l, t_{l+1}]$ and put

$$\begin{aligned} A &= a(t_l, X^{\mathrm{WP}}(t_l)) - a(t_l, X^{\mathrm{WPt}}(t_l)) \\ &\quad - a^{(0,1)}(t_l, X^{\mathrm{WPt}}(t_l)) \cdot (X^{\mathrm{WP}}(t_l) - X^{\mathrm{WPt}}(t_l)), \\ B &= \sigma(t_l, X^{\mathrm{WP}}(t_l)) - \sigma(t_l, X^{\mathrm{WPt}}(t_l)) \\ &\quad - \sigma^{(0,1)}(t_l, X^{\mathrm{WPt}}(t_l)) \cdot (X^{\mathrm{WP}}(t_l) - X^{\mathrm{WPt}}(t_l)) \end{aligned}$$

as well as

$$U = (X^{\mathrm{WP}}(t_l) - X^{\mathrm{aux}}(t_l)) \cdot (1 + a^{(0,1)}(t_l, X^{\mathrm{WPt}}(t_l)) \cdot (t - t_l))$$

and

$$\begin{aligned} V &= A \cdot (t - t_l) \\ &\quad + (\sigma^{(0,1)}(t_l, X^{\mathrm{WPt}}(t_l)) \cdot (X^{\mathrm{WP}}(t_l) - X^{\mathrm{aux}}(t_l)) + B) \cdot (W(t) - W(t_l)) \\ &\quad + (g_1(t_l, X^{\mathrm{WP}}(t_l)) - g_1(t_l, X^{\mathrm{WPt}}(t_l))) \cdot ((W(t) - W(t_l))^2 - (t - t_l)) \\ &\quad + (g_2(t_l, X^{\mathrm{WP}}(t_l)) - g_2(t_l, X^{\mathrm{WPt}}(t_l))) \cdot (W(t) - W(t_l)) \cdot (t - t_l) \\ &\quad + (g_3(t_l, X^{\mathrm{WP}}(t_l)) - g_3(t_l, X^{\mathrm{WPt}}(t_l))) \cdot (W(t) - W(t_l))^3 \\ &\quad + (g_4(t_l, X^{\mathrm{WP}}(t_l)) - g_4(t_l, X^{\mathrm{WPt}}(t_l))) \cdot (t - t_l)^2 \\ &\quad + (g_5(t_l, X^{\mathrm{WP}}(t_l)) - g_5(t_l, X^{\mathrm{WPt}}(t_l))) \cdot \int_{t_l}^{t} (W(s) - W(t_l))\, ds. \end{aligned}$$

By definition,

$$X^{\mathrm{WP}}(t) - X^{\mathrm{aux}}(t) = U + V.$$

Due to (A), we have

$$\begin{aligned} &|g_n(t_l, X^{\mathrm{WP}}(t_l)) - g_n(t_l, X^{\mathrm{WPt}}(t_l))| \\ &\qquad \le c \cdot (1 + |X^{\mathrm{WPt}}(t_l)|^2) \cdot |X^{\mathrm{WP}}(t_l) - X^{\mathrm{WPt}}(t_l)| \end{aligned}$$

for $n = 1, \ldots, 5$, and

$$|A| \le c \cdot |X^{\mathrm{WP}}(t_l) - X^{\mathrm{WPt}}(t_l)|^2, \qquad |B| \le c \cdot |X^{\mathrm{WP}}(t_l) - X^{\mathrm{WPt}}(t_l)|^2.$$

Put $q = \lceil p/2 \rceil$ and note that $p \le 2q \le 2p$. For $r = 1, \ldots, 2q$, we obtain

$$|U|^r \le (1 + c \cdot (t - t_l)) \cdot |X^{\mathrm{WP}}(t_l) - X^{\mathrm{aux}}(t_l)|^r$$

as well as

$$\begin{aligned} E(|V|^r | \mathcal{F}_{t_l}) &\le c \cdot |X^{\mathrm{WP}}(t_l) - X^{\mathrm{aux}}(t_l)|^r \cdot (t - t_l)^{r/2} \\ &\quad + c \cdot |X^{\mathrm{WP}}(t_l) - X^{\mathrm{WPt}}(t_l)|^{2r} \cdot (t - t_l)^{r/2} \\ &\quad + c \cdot (1 + |X^{\mathrm{WPt}}(t_l)|^{2r}) \cdot |X^{\mathrm{WP}}(t_l) - X^{\mathrm{WPt}}(t_l)|^r \cdot (t - t_l)^r. \end{aligned}$$

Moreover,

$$
\begin{aligned}
|E(V|\mathcal{F}_{t_l})| &= |A\cdot(t-t_l) + (g_4(t_l, X^{\mathrm{WP}}(t_l)) - g_4(t_l, X^{\mathrm{WPt}}(t_l)))\cdot(t-t_l)^2| \\
&\le c\cdot |X^{\mathrm{WP}}(t_l) - X^{\mathrm{WPt}}(t_l)|^2\cdot(t-t_l) \\
&\quad + c\cdot(1+|X^{\mathrm{WP}}(t_l)|^2)\cdot|X^{\mathrm{WP}}(t_l) - X^{\mathrm{WPt}}(t_l)|\cdot(t-t_l)^2.
\end{aligned}
$$

Hence,

$$
\begin{aligned}
&E(|X^{\mathrm{WP}}(t_l) - X^{\mathrm{aux}}(t_l)|^{2q}|\mathcal{F}_{t_l}) \\
&\quad = \sum_{r=0}^{2q}\binom{2q}{r}\cdot U^{2q-r}\cdot E(V^r|\mathcal{F}_{t_l}) \\
&\quad \le (1+c\cdot(t-t_l))\cdot|X^{\mathrm{WP}}(t_l) - X^{\mathrm{aux}}(t_l)|^{2q} \\
&\qquad + c\cdot\sum_{r=2}^{2q}\binom{2q}{r}\cdot|X^{\mathrm{WP}}(t_l) - X^{\mathrm{aux}}(t_l)|^{2q}\cdot(t-t_l)^{r/2} \\
&\qquad + c\cdot 2q\cdot|U|^{2q-1}\cdot(1+|X^{\mathrm{WP}}(t_l)|^2)\cdot|X^{\mathrm{WP}}(t_l) - X^{\mathrm{WPt}}(t_l)|\cdot(t-t_l)^2 \\
&\qquad + c\cdot\sum_{r=2}^{2q}\binom{2q}{r}\cdot|U|^{2q-r}\cdot(1+|X^{\mathrm{WP}}(t_l)|^{2r}) \\
&\qquad\qquad \times|X^{\mathrm{WP}}(t_l) - X^{\mathrm{WPt}}(t_l)|^r\cdot(t-t_l)^r \\
&\qquad + c\cdot\sum_{r=1}^{2q}\binom{2q}{r}\cdot|U|^{2q-r}\cdot|X^{\mathrm{WP}}(t_l) - X^{\mathrm{WPt}}(t_l)|^{2r}\cdot(t-t_l)^{r/2} \\
&\quad \le (1+c\cdot(t-t_l))\cdot|X^{\mathrm{WP}}(t_l) - X^{\mathrm{aux}}(t_l)|^{2q} \\
&\qquad + c\cdot\sum_{r=1}^{2q}\binom{2q}{r}\cdot|U|^{2q-r}\cdot(1+|X^{\mathrm{WP}}(t_l)|^{2r}) \\
&\qquad\qquad \times|X^{\mathrm{WP}}(t_l) - X^{\mathrm{WPt}}(t_l)|^r\cdot(t-t_l)^{1+r/2} \\
&\qquad + c\cdot\sum_{r=1}^{2q}\binom{2q}{r}\cdot|U|^{2q-r}\cdot|X^{\mathrm{WP}}(t_l) - X^{\mathrm{WPt}}(t_l)|^{2r}\cdot(t-t_l)^{\max(1,r/2)}.
\end{aligned}
$$

By Lemma 10, we get

$$
\begin{aligned}
&E(|U|^{2q-r}\cdot(1+|X^{\mathrm{WP}}(t_l)|^{2r})\cdot|X^{\mathrm{WP}}(t_l) - X^{\mathrm{WPt}}(t_l)|^r) \\
&\quad \le \|U\|_{2q}^{2q-r}\cdot(E((1+|X^{\mathrm{WP}}(t_l)|^{4q})\cdot|X^{\mathrm{WP}}(t_l) - X^{\mathrm{WPt}}(t_l)|^{2q}))^{r/(2q)} \\
&\quad \le \|U\|_{2q}^{2q-r}\cdot(1+\|X^{\mathrm{WP}}(t_l)\|_{8q}^{2r})\cdot\|X^{\mathrm{WP}}(t_l) - X^{\mathrm{WPt}}(t_l)\|_{4q}^{r} \\
&\quad \le c\cdot\|X^{\mathrm{WP}}(t_l) - X^{\mathrm{aux}}(t_l)\|_{2q}^{2q-r}\cdot\Delta_{\max}^{r}
\end{aligned}
$$

and similarly

$$
\begin{aligned}
&E(|U|^{2q-r} \cdot |X^{\mathrm{WP}}(t_l) - X^{\mathrm{WPt}}(t_l)|^{2r}) \\
&\quad \le \|U\|_{2q}^{2q-r} \cdot \|X^{\mathrm{WP}}(t_l) - X^{\mathrm{WPt}}(t_l)\|_{4q}^{2r} \\
&\quad \le c \cdot \|X^{\mathrm{WP}}(t_l) - X^{\mathrm{aux}}(t_l)\|_{2q}^{2q-r} \cdot \Delta_{\max}^{2r}.
\end{aligned}
$$

Thus,

$$
\begin{aligned}
&\|X^{\mathrm{WP}}(t) - X^{\mathrm{aux}}(t)\|_{2q}^{2q} \\
&\quad \le (1 + c \cdot (t - t_l)) \cdot \|X^{\mathrm{WP}}(t_l) - X^{\mathrm{aux}}(t_l)\|_{2q}^{2q} \\
&\qquad + c \cdot (t - t_l) \cdot \sum_{r=1}^{2q} \binom{2q}{r} \|X^{\mathrm{WP}}(t_l) - X^{\mathrm{aux}}(t_l)\|_{2q}^{2q-r} \cdot \Delta_{\max}^{3r/2} \\
&\quad \le (1 + c \cdot (t - t_l)) \cdot \|X^{\mathrm{WP}}(t_l) - X^{\mathrm{aux}}(t_l)\|_{2q}^{2q} \\
&\qquad + c \cdot (t - t_l) \cdot (\|X^{\mathrm{WP}}(t_l) - X^{\mathrm{aux}}(t_l)\|_{2q} + \Delta_{\max}^{3/2})^{2q} \\
&\quad \le (1 + c \cdot (t - t_l)) \cdot \|X^{\mathrm{WP}}(t_l) - X^{\mathrm{aux}}(t_l)\|_{2q}^{2q} + c \cdot (t - t_l) \cdot \Delta_{\max}^{3q}.
\end{aligned}
$$

Now apply Gronwall's lemma to complete the proof. □

**Acknowledgments.** I thank the referees for a number of constructive suggestions that helped to improve the presentation of the paper. I am also grateful to Klaus Ritter and Norbert Hofmann for very stimulating discussions.

Fakultät für Mathematik
Institut für Mathematische Stochastik
Otto-von-Guericke-Universität Magdeburg
Postfach 4120, 39016 Magdeburg
Germany
e-mail: gronbach@mail.math.uni-magdeburg.de